\documentclass[a4paper]{amsart}
\usepackage{pictex}
\usepackage{amscd} 
\usepackage{ifthen}
\usepackage{graphics}
\usepackage{amsmath}
\usepackage{latexsym}
\usepackage{amssymb}

\begin{document}

\newcommand{\sH}{\operatorname{\mathcal H}\nolimits}
\newcommand{\Tilt}{\operatorname{Tilt}\nolimits}
\newcommand{\soc}{\operatorname{soc}\nolimits}
\newcommand{\Tr}{\operatorname{Tr}\nolimits}
\newcommand{\ind}{\operatorname{ind}\nolimits}
\newcommand{\Supp}{\operatorname{Supp}\nolimits}
\newcommand{\extto}{\xrightarrow}
\renewcommand{\mod}{\operatorname{mod}\nolimits}
\newcommand{\Mod}{\operatorname{Mod}\nolimits}
\newcommand{\Sub}{\operatorname{Sub}\nolimits}
\newcommand{\add}{\operatorname{add}\nolimits}
\newcommand{\Hom}{\operatorname{Hom}\nolimits}
\newcommand{\Rad}{\operatorname{Rad}\nolimits}
\newcommand{\End}{\operatorname{End}\nolimits}
\renewcommand{\Im}{\operatorname{Im}\nolimits}
\newcommand{\Ker}{\operatorname{Ker}\nolimits}
\newcommand{\Coker}{\operatorname{Coker}\nolimits}
\newcommand{\rad}{\operatorname{rad}\nolimits}
\newcommand{\Ext}{\operatorname{Ext}\nolimits}
\newcommand{\op}{{\operatorname{op}}}
\newcommand{\Ab}{\operatorname{Ab}\nolimits}
\newcommand{\id}{\operatorname{id}\nolimits}
\newcommand{\pd}{\operatorname{pd}\nolimits}
\renewcommand{\P}{\operatorname{\mathcal P}\nolimits}
\newcommand{\A}{\operatorname{\mathcal A}\nolimits}
\newcommand{\C}{\operatorname{\mathcal C}\nolimits}
\newcommand{\D}{\operatorname{\mathcal D}\nolimits}
\newcommand{\E}{\operatorname{\mathcal E}\nolimits}
\newcommand{\G}{\Gamma}
\renewcommand{\L}{\Lambda}
\newcommand{\N}{\operatorname{\mathbb N}\nolimits}
\renewcommand{\S}{\operatorname{\mathcal S}\nolimits}
\newcommand{\T}{\operatorname{\mathcal T}\nolimits}
\newcommand{\Z}{\operatorname{\mathbb Z}\nolimits}
\newtheorem{lem}{Lemma}[section]
\newtheorem{lemma}[lem]{Lemma}
\newtheorem{prop}[lem]{Proposition}
\newtheorem{cor}[lem]{Corollary}
\newtheorem{thm}[lem]{Theorem}
\newtheorem{theorem}[lem]{Theorem}
\newtheorem{example}[lem]{Example}
\newtheorem{rem}[lem]{Remark}
\newtheorem{remark}[lem]{Remark}
\newtheorem{defin}[lem]{Definition}
\newtheorem{definition}[lem]{Definition}
\newtheorem{ex}[lem]{Example}
\newtheorem{conj}[lem]{Conjecture}
\newtheorem{problem}[lem]{Problem}

\title[Tilting and clusters]{Tilting theory and cluster combinatorics}

\author[Buan]{Aslak Bakke Buan}
\address{Institutt for matematiske fag\\
Norges teknisk-naturvitenskapelige universitet\\
N-7491 Trondheim\\
Norway}
\email{aslakb@math.ntnu.no}

\author[Marsh]{Robert Marsh}
\address{Department of Mathematics \\
University of Leicester \\
University Road \\
Leicester LE1 7RH \\
England
}
\email{rjm25@mcs.le.ac.uk}

\author[Reineke]{Markus Reineke}
\address{BU Wuppertal \\
Gausstrasse 20 \\
D-42097 Wuppertal
Germany}
\email{reineke@math.uni-wuppertal.de}

\author[Reiten]{Idun Reiten}
\address{Institutt for matematiske fag\\
Norges teknisk-naturvitenskapelige universitet\\
N-7491 Trondheim\\
Norway}
\email{idunr@math.ntnu.no}

\author[Todorov]{Gordana Todorov}
\address{Department of mathematics\\
Northeastern University\\
Boston\\
MA 02115\\
USA }
\email{todorov@neu.edu}

\keywords{APR tilting theory, tilting module, Hom-configuration,
Ext-configuration, complement, approximation theory,
self-injective algebra, cluster algebra.}
\subjclass[2000]{Primary: 16G20, 16G70; Secondary 16S99, 17B99}

\begin{abstract}
We introduce a new category $\C$, which we call the {\em cluster category},
obtained as a quotient of the bounded derived category $\D$ of the module
category of a finite-dimensional hereditary algebra $H$ over a field.
We show that, in the simply-laced
Dynkin case, $\C$ can be regarded as a natural model for the combinatorics of
the corresponding Fomin--Zelevinsky cluster algebra. In this model,
the tilting modules correspond to the clusters of Fomin--Zelevinsky. 
Using approximation theory, we investigate the tilting theory of $\C$,
showing that it is more regular than that of the module category itself, and
demonstrating an interesting link with the classification of self-injective
algebras of finite representation type. This investigation also enables us to
conjecture a generalisation of APR-tilting. 
\end{abstract}

\date{4th February 2004}

\thanks{Aslak Bakke Buan was supported by a grant
from the Norwegian Research Council.
Robert Marsh was supported by a Leverhulme Fellowship,
study leave from the University of Leicester,
EPSRC grant GR/R17547/01, NSF funding through a research grant of
A.~Zelevinsky and support from the Norwegian Research Council.}

\maketitle

\section*{Introduction}
In this paper, we introduce a new category, which we call the {\em cluster
category}, associated with any finite dimensional hereditary algebra $H$
over a field $k$. This is defined as the quotient $\C$ of the
bounded derived category $\D$ of finitely generated modules over
$H$ by the functor $F=\tau^{-1}[1]$, where $\tau$ denotes the AR-translation
and $[1]$ denotes the shift functor. The category $\C$ is
triangulated, by a result of Keller~\cite{k}, and we show that it is also a
Krull-Schmidt category. Our main aims are to show how this category can be
used to study the tilting theory of $H$ (and related algebras) and to show
that it can be used as a model for the combinatorics of an associated
Fomin--Zelevinsky~\cite{fz1} cluster algebra.

{\em $\Hom$-configurations} are certain collections of non-isomorphic
indecomposable objects in $\D$, and were considered in~\cite{rie1} in
connection with the classification of self-injective algebras of finite
representation type. We formulate analogous conditions using $\Ext^1$ instead
of $\Hom$, and call the resulting collections {\em $\Ext$-configurations}.
We show that they exhibit a behaviour similar to that of $\Hom$-configurations.
In particular, they are invariant under the functor $F$
(compare~\cite{blr}, where it is shown that $\Hom$-configurations exhibit a
similar kind of invariance in the Dynkin case). As a consequence, we can
show that they are in 1--1 correspondence with basic tilting objects in $\C$.
By showing that a basic tilting object in $\C$ is induced by a basic tilting
module over some hereditary algebra derived equivalent to $H$, we prove that
$\Ext$-configurations, like $\Hom$-configurations in the Dynkin case,
are induced by basic tilting modules.

The category $\C$ provides an interesting ``extension'' of the module
category of $H$.
It is known that any almost complete basic tilting module $\overline{T}$ over
$H$ can be completed to a basic tilting module in at most two different
ways~\cite{rs1,u1} and in exactly two different ways if and only if
$\overline{T}$ is sincere~\cite{hu1}.
However, in the extended category $\C$,
the behaviour is more regular: an almost complete basic tilting object always
has exactly two complements. We show further that, given one complement $M$ to
an almost complete basic tilting object $\overline{T}$,
the other can be constructed using approximation theory from~\cite{as}.
Indeed, we show that there is a triangle
$$M^{\ast}\rightarrow B\rightarrow M\rightarrow M^{\ast}[1]$$
in $\C$, where $B\rightarrow M$ is a minimal right
$\add\overline{T}$-approximation of $M$ in $\C$ and $M^{\ast}$ is
the other complement to $\overline{T}$. Dually, there is a triangle
$$M\rightarrow B'\rightarrow M^{\ast}\rightarrow M[1]$$
in $\C$. In fact, we are able to show that two indecomposable objects
$M$ and $M^{\ast}$ form such an exchange pair if and only if
$$\dim_{\End(M)}\Ext^1_{\C}(M,M^{\ast}) = 1 = 
\dim_{\End(M^{\ast})}\Ext^1_{\C}(M^{\ast},M).$$

The above results have some interesting interpretations in the Dynkin
case in terms of cluster algebras, which were defined by Fomin and
Zelevinsky~\cite{fz1}. These algebras were defined so that the cluster
structure (when quantised) should encode multiplicative properties of the dual
canonical basis of the quantised enveloping algebra of a semisimple Lie
algebra over $\mathbb{C}$, and that it should model the (classical and
quantised) coordinate rings of varieties associated to algebraic groups
(now shown in several cases --- see~\cite{fz3},~\cite{s}.),
with particular relevance to total positivity properties; there have already
been many applications to other areas as
well~\cite{cfz,fz3,fz4,fz5,gsv,mrz,p}.

The definition is as follows.
Let $\mathbb{F}=\mathbb{Q}(u_1,u_2,\ldots ,u_n)$ be the field
of rational functions in indeterminates $u_1,u_2,\ldots u_n$. Let
$\mathbf{x}\subseteq \mathbb{F}$ be a transcendence basis over $\mathbb{Q}$,
and let
$B=(b_{xy})_{x,y\in \mathbf{x}}$ be an $n\times n$ sign-skew-symmetric
integer matrix with rows and columns indexed by $\mathbf{x}$. In other words,
we suppose that for all $x,y\in \mathbf{x}$, $b_{xy}=0$ if and only if
$b_{yx}=0$, that $b_{xy}>0$ if and only if $b_{yx}<0$, and that $b_{xx}=0$.
Such a pair $(\mathbf{x},B)$ is called a {\em seed}. Fomin and
Zelevinsky~\cite{fz1},~\cite{fz2} have defined a certain subring
$\mathcal{A}(\mathbf{x},B)$ of $\mathbb{F}$ associated to the seed
$(\mathbf{x},B)$, known as a {\em cluster algebra}. Given such a seed,
and an element $z\in \mathbf{x}$, define a new element $z'\in\mathbb{F}$ via
the {\em binary exchange relation}:
\begin{equation}\label{exchangerelation}
zz'=\prod_{x\in\mathbf{x},b_{xz}>0}x^{b_{xz}}+
\prod_{x\in\mathbf{x},b_{xz}<0}x^{-b_{xz}}.
\end{equation}
In such circumstances, we say that $z,z'$ form an {\em exchange pair}.
Let $\mathbf{x}'=\mathbf{x}\cup \{z'\}\setminus \{z\}$, a new
transcendence basis of $\mathbb{F}$. Let $B'$ be the
{\em mutation} of the matrix $B$ in direction $z$
(as defined in~\cite{fz1}). Then
$$b'_{xy}=\left\{
\begin{array}{ll}
-b_{xy} & \mbox{if\ } x=z \mbox{\ or\ }y=z, \\
b_{xy}+\frac{1}{2}(|b_{xz}|b_{zy}+b_{xz}|b_{zy}|) & \mbox{otherwise.}
\end{array}
\right. $$
The row and column labelled $z$ in $B$ are relabelled $z'$ in
$B'$. The pair $(\mathbf{x}',B')$ is called the {\em mutation} of the seed
$\mathbf{x}$ in direction
$z$. Let $\mathcal{S}$ be the set of seeds obtained by iterated
mutation of $(\mathbf{x},B)$. Then the set of
{\em cluster variables} is, by definition, the union $\chi$ of the
transcendence bases appearing in the seeds in $\mathcal{S}$, and the
cluster algebra $\mathcal{A}(\mathbf{x},B)$ is the subring of $\mathbb{F}$
generated by $\chi$. Up to isomorphism of cluster algebras, it does not
depend on the initial choice $\mathbf{x}$ of transcendence basis,
so can be denoted $\mathcal{A}_B$.
In general, coefficients appear in the relation~\eqref{exchangerelation},
but here we take all of these coefficients to be $1$ as this is enough
to describe the connections with representation theory that we consider.

If its matrix is skew-symmetric, a seed $(\mathbf{x},B)$ determines a
quiver with vertices corresponding to its rows and columns, and $b_{ij}$
arrows from vertex $i$ to vertex $j$ whenever $b_{ij}>0$.
If $\chi$ is finite, the cluster algebra $\mathcal{A}_B$ is said to be of
finite type. In~\cite{fz2}, it is shown that, up to isomorphism,
the cluster algebras of finite type can be classified by the Dynkin diagrams;
they are precisely those for which there exists a seed whose corresponding
quiver is of Dynkin type. In this case, Fomin and Zelevinsky associate
a nonnegative integer, known as the {\em compatibility degree}, to each pair
of cluster variables (see Section~\ref{connections}).
Two variables are said to be {\em compatible}
provided that their compatibility degree is zero, and clusters are maximal
compatible subsets of $\chi$.

Suppose that $H$ is the path algebra of a simply-laced Dynkin quiver of type
$\Delta$. We show that the indecomposable objects in $\C$ are in 1--1
correspondence with the cluster variables in a cluster algebra $\A$ of type
$\Delta$. Using results from~\cite{mrz} we show that, for the two
indecomposable objects $X,Y$ in $\C$, $\dim \Ext^1_{\C}(X,Y)$ is equal to
the compatibility degree of the corresponding cluster variables.

The advantage of our approach here is that it allows us to give a
direct interpretation of all clusters in terms of tilting objects:
it follows from the above that the clusters of $\A$ are in 1--1 correspondence
with the basic tilting objects in $\C$. We develop this relationship
further: the existence of exactly two complements for
any almost complete basic tilting object in $\C$ then corresponds to the
fact that for any almost complete cluster there are exactly two ways to
complete it to a cluster (by adding a cluster variable). A consequence of
our result above is a new proof of the result~\cite[3.5,4.4]{fz2} that
two cluster variables form an exchange pair
(i.e.\ appear in an exchange relation --- see
equation~\eqref{exchangerelation}) if and only if their compatibility
degree is $1$. We conjecture that in this case the middle term $B$ in the
triangle above is the direct sum of the indecomposable objects corresponding
to the cluster variables appearing in one term of the exchange
relation~\cite[1.1]{fz2}, with the middle term $B'$ of the
dual triangle corresponding to the other term
(see Conjecture~\ref{exchangeconjecture}), suggesting that it might be
possible to construct the cluster algebra directly from $\C$.
Finally, we are able to use the new perspective on tilting theory
afforded by cluster algebras and the cluster category
to conjecture a generalisation of APR-tilting (see~\cite{apr}).

P.~Caldero, F.~Chapoton and R.~Schiffler~\cite{ccs} have recently
associated a category to the cluster algebra of type $A_n$, giving
a definition via the combinatorics of the corresponding cluster algebra.
They have shown that this category is equivalent to the cluster category
$\C$ we have associated to a Dynkin quiver of type $A_n$.
Their approach enables them to generalise the
denominator theorem of Fomin and Zelevinsky~\cite[1.9]{fz2}
to an arbitrary cluster.
Instead, in our approach we consider a more general situation
(an arbitrary finite dimensional hereditary algebra), and the connections
with tilting theory and configurations of modules in the derived
category. We develop links with cluster combinatorics for all simply-laced
Dynkin cases in a uniform way.

\section{Cluster categories} \label{clustercategories}

In this section we introduce what we call the cluster category of a finite 
dimensional hereditary algebra, and discuss some of its elementary properties.

Let $H$ be a finite dimensional hereditary algebra over a field $k$, and
denote by $\D = D^b(H)$ the bounded derived category
of finitely generated $H$-modules with shift functor $[1]$. For any category
$\E$, we will denote by $\ind\E$ the subcategory of
isomorphism classes of indecomposable objects in $\E$; depending on the
context we shall also use the same notation to denote the set of isomorphism
classes of indecomposable objects in $\E$.

Let $G \colon \D \to \D$ be a triangle functor, which we also assume 
satisfies the following properties; see~\cite{k}.
\begin{itemize}
\item[(g1)]{For each $U$ in $\ind H$, only a finite number 
of objects $G^n U$, where $n \in \Z$, lie in $\ind H$.}
\item[(g2)]{There is some $N \in \N$ such that 
$\{U[n] \mid U \in \ind H, n \in [-N,N] \}$ contains 
a system of representatives of the orbits of $G$ on $\ind \D$.}
\end{itemize}

We denote by $\D/ G$ the corresponding factor category. The objects are by 
definition the $G$-orbits of objects in $\D$, and the
morphisms are given by 
$$\Hom_{\D/G}(\widetilde{X},\widetilde{Y}) = \coprod_{i \in \Z} \Hom_{\D}(G^iX,Y).$$
Here $X$ and $Y$ are objects in $\D$, and $\widetilde{X}$ and $\widetilde{Y}$
are the corresponding objects in $\D/G$ (although we shall often write
such objects simply as $X$ and $Y$).
Note that it follows from our assumptions on $G$ that
$\Hom_{\D}(G^i X,Y) \neq 0$ for only a finite number of values
of $i$. It is known from~\cite{k} that $\D/G$ is a triangulated
category and that 
the natural functor $\pi \colon \D \to \D/G$ is a triangle functor.
The shift in $\D/G$ is induced by the shift in $\D$, and is also denoted
by $[1]$. In both cases we write as usual $\Hom(U,V[1]) = \Ext^1(U,V)$.
We then have 
$$\Ext^1_{\D/G}(\widetilde{X},\widetilde{Y}) = \coprod_{i \in \Z} 
\Ext^1_{\D}(G^i X,Y),$$
where $X,Y$ are objects in $\D$ and $\widetilde{X},\widetilde{Y}$ are the
corresponding objects in $\D/G$.
Note that since there are only finitely many values of $i$ such that
$\Hom_{\D}(G^i X,Y)$ is not zero, there
are also only finitely many values of $i$ such that
$\Ext^1_{\D}(G^iX,Y)$ is not zero, for $X,Y$ in $\D$.
We remark that the quotient $D^b(H)/[2]$ was
considered in~\cite{h1}; however, this quotient has quite
different properties and is not closely linked with cluster algebras.

While several properties hold for arbitrary functors $G$ satisfying
(g1) and (g2), we shall
mainly be concerned with the special choice of functor $F = \tau^{-1}[1]$,
where $\tau$ is the AR-translation in $\D$
(which is induced by $D\Tr$ on non-projective indecomposable 
objects in $\ind H$, and where $\tau(P) = I[-1]$ when $P$
is indecomposable projective and $I$ denotes the indecomposable
injective with $\soc I \simeq P/\underline{r}P$).

We shall see various reasons why the factor category $\D / F$ 
is especially nice. Because of the applications to cluster 
theory we call it the \emph{cluster category} of $H$,
and we denote it by $\C$. 

If we are in the setting with $H$ of finite representation type and
$k$ an algebraically closed field, then $\D$ (and thus $\C$) only depends 
on the underlying graph $\Delta$ of the quiver of $H$, and we write $\C =
\C(\Delta)$. Then $\Delta$ is a simply-laced Dynkin 
diagram. For this case we give a combinatorial construction of
$\ind\C$. We recall the theory of translation quivers
from~\cite{rie2}. If $\Gamma=(\Gamma_0,\Gamma_1)$ is any quiver, with
vertices $\Gamma_0$ and arrows $\Gamma_1$, we recall that, if $x\in \Gamma_0$,
then $x^+$ is the set of the end-points of arrows which start at $x$, while
$x^-$ denotes the set of starting points of arrows which end at $x$.
A {\em stable translation quiver} is
a quiver $\Gamma$, without any loops or multiple edges, together with a
bijection $\tau:\Gamma_0\rightarrow \Gamma_0$ (known as the {\em translation})
such that, for all $x\in\Gamma_0$, $x^-=\tau(x)^+$. A morphism
of stable translation quivers is defined to be a quiver morphism which
commutes with translation.

If $\Gamma$ is a stable translation quiver, and $a:x\rightarrow y$ is an arrow
of $\Gamma$, then there is a unique arrow $\sigma(a):\tau(y)\rightarrow x$.
The rule $a\mapsto \sigma(a)$ defines a bijection from $\Gamma_1$
to $\Gamma_1$, known as the {\em polarisation}.
The {\em mesh category} associated to $\Gamma$ has objects indexed by
the vertices of $\Gamma$, and morphisms generated by the arrows of
$\Gamma$, subject to the mesh relations (for all vertices $y$ of
$\Gamma$):
$$
\sum_{a:x\rightarrow y} \sigma(a)a=0.
$$
If $\mathcal{E}$ is a Krull-Schmidt category with almost split sequences,
we shall denote its
AR-quiver by $\Gamma(\mathcal{E})$ (see~\cite{rin}).

Let $Q$ be the quiver of $H$ and
let $\mathbb{Z}Q$ be the stable translation quiver associated to
$Q$ (see~\cite{rie2}).
The vertices of $\mathbb{Z}Q$ are labelled by pairs $(n,i)$ with
$n$ in $\mathbb{Z}$ and $i$ a vertex of $Q$. Whenever there is an arrow in
$Q$ from $i$ to $j$ there is an arrow from $(n,i)$ to $(n,j)$, and an
arrow from $(n,j)$ to $(n+1,i)$, and these are all the arrows
in $\mathbb{Z}Q$. A translation $\tau$ is defined on
$\mathbb{Z}Q$, just taking $(n,i)$ to $(n-1,i)$. In this way
$\mathbb{Z}Q$ is a stable translation quiver.
We denote the corresponding mesh category by $k(\mathbb{Z}Q)$.
We have the following:

\begin{prop} (Happel~\cite[5.6]{h2}) \label{happelequivalence}
Let $Q$ be any quiver of Dynkin type. Then the mesh category
$k(\mathbb{Z}Q)$ is equivalent to $\ind\D$.
\end{prop}

It follows that (as a stable translation quiver), $\mathbb{Z}Q$
depends only on the underlying Dynkin diagram $\Delta$, and not on $Q$.
We therefore denote it $\mathbb{Z}\Delta$, and denote the corresponding
mesh category by $k(\mathbb{Z}\Delta)$.
The AR-quiver of $\D$ is $\Gamma(\D)=\mathbb{Z}\Delta$.

We recall that $F=\tau^{-1}[1]$ is an autoequivalence of $\D$, and
therefore permutes the indecomposable objects, inducing a graph automorphism
$\varphi$ (via Proposition~\ref{happelequivalence}) of $\mathbb{Z}\Delta$.
We note that the graph automorphisms induced by $\tau^{-1}$ and $[1]$
are independent of the orientation $Q$, so $\varphi$ is
independent of $Q$. Since $F$ commutes with $\tau$ on $\D$, $\varphi$ is an
automorphism of stable translation quivers. It follows that the quotient
graph $\mathbb{Z}\Delta/\varphi$ is also a stable translation quiver, and we
can form the corresponding mesh category; this is equivalent to the category
$\ind\C(\Delta)$ defined above.
The natural epimorphism of stable translation quivers
$\pi:\mathbb{Z}\Delta\rightarrow \mathbb{Z}\Delta/\varphi$, taking
the vertex $v$ of $\mathbb{Z}\Delta$ to its $\varphi$-orbit $\pi(v)$,
induces the functor $\pi$ above.

\subsection*{Example}
\begin{figure}[htbp]
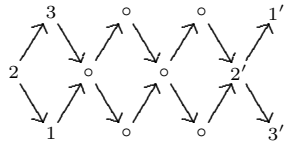


\beginpicture

\setcoordinatesystem units <0.5cm,0.4cm>             
\setplotarea x from -8 to 10, y from -9 to -3       

\scriptsize{

\put{$1$}[c] at 2 -8
\put{$\circ$}[c] at 4 -8
\put{$\circ$}[c] at 6 -8
\put{$3'$}[c] at 8 -8
\put{$2$}[c] at 1 -6
\put{$\circ$}[c] at 3 -6
\put{$\circ$}[c] at 5 -6
\put{$2'$}[c] at 7 -6
\put{$3$}[c] at 2 -4
\put{$\circ$}[c] at 4 -4
\put{$\circ$}[c] at 6 -4
\put{$1'$}[c] at 8 -4
}
\setlinear \plot  2 -7.65   2.6 -6.4  / %
\setlinear \plot  2.6 -6.4 2.7 -6.7  / %
\setlinear \plot  2.6 -6.4 2.3 -6.5  / %

\setlinear \plot  4 -7.65   4.6 -6.4  / %
\setlinear \plot  4.6 -6.4 4.7 -6.7  / %
\setlinear \plot  4.6 -6.4 4.3 -6.5  / %

\setlinear \plot  6 -7.65   6.6 -6.4  / %
\setlinear \plot  6.6 -6.4 6.7 -6.7  / %
\setlinear \plot  6.6 -6.4 6.3 -6.5  / %

\setlinear \plot  1.6 -7.65 1.0 -6.4  / %
\setlinear \plot  1.6 -7.65 1.3 -7.55  / %
\setlinear \plot  1.6 -7.65 1.7 -7.35 / %

\setlinear \plot  3.6 -7.65 3.0 -6.4  / %
\setlinear \plot  3.6 -7.65 3.3 -7.55 / %
\setlinear \plot  3.6 -7.65 3.7 -7.35  / %

\setlinear \plot  5.6 -7.65 5.0 -6.4  / %
\setlinear \plot  5.6 -7.65 5.3 -7.55  / %
\setlinear \plot  5.6 -7.65 5.7 -7.35  / %

\setlinear \plot  7.6 -7.65 7.0 -6.4  / %
\setlinear \plot  7.6 -7.65 7.3 -7.55  / %
\setlinear \plot  7.6 -7.65 7.7 -7.35  / %

\setlinear \plot  1   -5.6 1.6 -4.4  / %
\setlinear \plot  1.6 -4.4 1.7 -4.7  / %
\setlinear \plot  1.6 -4.4 1.3 -4.5  / %

\setlinear \plot  3   -5.6 3.6 -4.4  / %
\setlinear \plot  3.6 -4.4 3.7 -4.7  / %
\setlinear \plot  3.6 -4.4 3.3 -4.5  / %

\setlinear \plot  5   -5.6 5.6 -4.4  / %
\setlinear \plot  5.6 -4.4 5.7 -4.7  / %
\setlinear \plot  5.6 -4.4 5.3 -4.5  / %

\setlinear \plot  7   -5.6 7.6 -4.4  / %
\setlinear \plot  7.6 -4.4 7.7 -4.7  / %
\setlinear \plot  7.6 -4.4 7.3 -4.5  / %

\setlinear \plot  2.6 -5.65 2.0 -4.4  / %
\setlinear \plot  2.6 -5.65 2.3 -5.55  / %
\setlinear \plot  2.6 -5.65 2.7 -5.35 / %

\setlinear \plot  4.6 -5.65 4.0 -4.4  / %
\setlinear \plot  4.6 -5.65 4.3 -5.55  / %
\setlinear \plot  4.6 -5.65 4.7 -5.35 / %

\setlinear \plot  6.6 -5.65 6.0 -4.4  / %
\setlinear \plot  6.6 -5.65 6.3 -5.55  / %
\setlinear \plot  6.6 -5.65 6.7 -5.35 / %

\endpicture
\caption{The AR-quiver of $\C$ in type $A_3$}
\label{Cquiver}
\end{figure}

In Figure~\ref{Cquiver}, we show the AR-quiver of $\C$ in
type $A_3$.
The objects $1$,$2$ and $3$ are identified with $1'$, $2'$ and $3'$
(so that, in some sense, the quotient is a M\"{o}bius strip).

We are mostly interested in the factor $\C = \D / F$, where
$F = \tau^{-1}[1]$. The next properties, however, we state and prove
in a more general setting. 

\begin{prop} \label{krullschmidt}
Let $\D = D^b(H)$ for a finite dimensional hereditary $k$-algebra
$H$, and let $G \colon \D \to \D$ be a triangle
functor satisfying (g1) and (g2). Then the triangulated category $\D/G$ 
is a Krull-Schmidt category.
\end{prop}

\begin{proof}
Let $\widetilde{X}$ be in $\D/G$ induced by $X$ in $\D$.
We know that $X = X_1 \coprod \cdots \coprod X_n$ in $\D$, where each $X_i$ 
is indecomposable, with local endomorphism ring. Since the functor
$\pi \colon \D \to \D/G$ commutes with finite direct sums, we
have $\widetilde{X} = \widetilde{X_1} \coprod \cdots \coprod \widetilde{X_n}$.

We then claim that $\End_{\D/G}(\widetilde{X_i})$ is local for 
each $i$. So let $Y$ be in $\ind \D$. By definition,
$\Hom_{\D/G}(\widetilde{Y},\widetilde{Y}) = 
\coprod_{i \in \Z} \Hom_{\D}(G^i Y,Y)$. It is easy to see that
$$\rad(Y,Y) \coprod(\coprod_{i \neq 0} \Hom_{\D}(G^i Y,Y))$$
is a unique maximal ideal in $\Hom_{\D/G}(\widetilde{Y},\widetilde{Y})$, which
is hence a local ring. Thus, $\D/G$ is a Krull-Schmidt category.
\end{proof}

We remark that triangles in $\D/G$ are not necessarily induced by those
in $\D$. However, we have the following:

\begin{prop}
Let $\D = D^b(H)$ for a finite dimensional hereditary $k$-algebra
$H$, and let $G \colon  \D \to \D$ be a triangle 
functor satisfying (g1) and (g2).
Then $\D/G$ has almost split triangles induced by those in $\D$, 
and the AR-quiver is $\G(\D)/\varphi(G)$, where $\varphi$ is
the graph automorphism induced by $G$.
\end{prop}

\begin{proof}
Let $\widetilde{X}$ be an indecomposable object in $\D/G$, induced
by $X$ in $\D$. Let $$\tau X \overset{f}{\rightarrow} E \overset{g}{\rightarrow}
X \overset{s}{\rightarrow} \tau X[1]$$
be an almost split triangle in $\D$. Since $\pi \colon \D \to \D /G$ is a triangle functor,
there is the induced triangle
\begin{equation} \label{inducedtriangle}
\widetilde{\tau X} \overset{\widetilde{f}}{\rightarrow} \widetilde{E} \overset{\widetilde{g}}{\rightarrow}
\widetilde{X} \overset{\widetilde{s}}{\rightarrow} \widetilde{\tau X}[1]
\end{equation}
in $\D/G$. Since $s \neq 0$, we clearly have $\widetilde{s} \neq 0$. Let $\widetilde{Z}$ be
in $\ind \D/G$, induced by $Z$ in $\ind \D$, with $\widetilde{Z} \not \simeq \widetilde{X}$, 
and let
$\widetilde{h} \colon \widetilde{Z} \to \widetilde{X}$ be nonzero. Then $\widetilde{h} = \coprod \; h_i$, with
$h_i \in \Hom_{\D}(G^iZ,X)$. Since $\widetilde{Z} \not \simeq \widetilde{X}$, we have $G^i Z \not \simeq X$ 
for all $i$ and hence there is some $t_i \colon G^i Z \to E$ such that $gt_i=h_i$.
Let $\widetilde{t} = \coprod \; t_i$. Then we have $\widetilde{g} \widetilde{t} = \widetilde{h}$,
and hence $\widetilde{g}$ is right almost split. Similarly $\widetilde{f}$ is 
a left almost split map, and hence~(\ref{inducedtriangle}) is an almost
split triangle, and the translation $\widetilde{\tau}$ in $\D /G$ is given by 
$\widetilde{\tau} X= \widetilde{\tau X}$.
Hence it also follows that $\G(\D)/\varphi$ is the AR-quiver for $\D/G$.
\end{proof}

Let $D = \Hom_k(\ ,k)$.
It is also useful to note that the Serre duality formula 
$D \Ext_{\D}^1(A,B) \simeq \Hom_{\D}(A,\tau B)$, valid in $D^b(H)$, induces an analogous
formula for $\D/G$.

\begin{prop}
Let the notation and assumptions be as above. Then for $\widetilde{X}$ and $\widetilde{Y}$
in $\D^b(H)/G$ we have the Serre duality formula:
$$D \Ext^1_{\C}(\widetilde{X},\widetilde{Y}) \simeq \Hom_{\C}(\widetilde{Y},\widetilde{\tau} \widetilde{X})$$
functorial in both $\widetilde{X}$ and $\widetilde{Y}$.
\end{prop}

\begin{proof}
We have $$\Ext^1_{\D/G}(\widetilde{X},\widetilde{Y}) = \coprod_i \Hom_{\D}(G^i X, Y[1]) =
\coprod_i \Ext^1_{\D}(G^i X,Y)$$ and $\Hom_{\D/G}(\widetilde{Y},\widetilde{\tau} \widetilde{X}) =
\coprod_i \Hom_{\D}(G^iY, \tau X)$. We then apply the corresponding formula
for $D^b(H)$. 
\end{proof}

We end this section with some properties of $\C$.
Let $\S = \ind(\mod H \vee H[1])$, i.e.\ the set consisting of the
indecomposable $H$-modules, together with the objects $P[1]$, where
$P$ is an indecomposable projective $H$-module. Then it can be seen that
$\S$ is a fundamental domain for the action of $F$ on $\ind \D$,
containing exactly one representative from each $F$-orbit on $\ind \D$.
We recall that there is an oriented graph structure on $\ind \D$,
with an arrow from object $X$ to object $Y$ if there is a non-zero
map from $X$ to $Y$.

\begin{prop} \label{fewnonzero}
Let $X$ and $Y$ be objects in $\S$.
\begin{itemize}
\item[(a)]{We have $\Hom_{\D}(F^iX,Y)=0$ for all $i \neq -1,0$.}
\item[(b)]{If $X$ or $Y$ does not lie on an oriented cycle in $\D$,
then $\Hom_{\D}(F^iX,Y)\not= 0$ for at most one value of $i$.}
\end{itemize}
\end{prop}

\begin{proof}
(a) We have $\Hom_{\D}(F^iX,Y) = \Hom_{\D}(\tau^{-i}X[i],Y)$. For $i \geq 1$,
we clearly have $\Hom_{\D}(\tau^{-i}X[i],Y)= 0$. This is obvious for $i>1$, 
and for $i=1$ we only have to consider the case $Y = P[1]$ for
$P$ an indecomposable projective $H$-module. In that case we have
$\Hom_{\D}(\tau^{-1}X,P)$, which must be 0.
For $i \leq -2$ we have that
$$\Hom_{\D}(\tau^{-i}X[i],Y) = \Ext_{\D}^{-i}(\tau^{-i}X,Y)= 0.$$
\noindent (b)
We have \sloppy
$\Hom_{\D}(F^{-1}X,Y) = \Hom_{\D}(\tau X[-1],Y) = \Ext_{\D}^1(\tau X,Y) \simeq 
D \Hom_{\D}(Y, \tau^2 X) \simeq D\Hom_{\D}(\tau^{-2}Y,X)$. 
If $\Hom_{\D}(X,Y) \neq 0$ and $\Hom_{\D}(\tau^{-2}Y,X) \neq 0$ then it is clear
that $X$ and $Y$ lie on a cycle. 
\end{proof}

\begin{prop} \label{indecomposableobjects}
The indecomposable objects in $\C$ are precisely those of the form
$\widetilde{X}$ for $X$ an object in $\S$.
\end{prop}

\begin{proof}
It follows from Proposition~\ref{krullschmidt} and its proof that
the objects $\widetilde{X}$ for $X$ an object in $S$ are indecomposable
objects in $\C$. Using the definition of morphisms in $\C$ it is easy to
see, using Proposition~\ref{fewnonzero}(a), that if $X,Y\in S$ are such that
$\widetilde{X}\simeq \widetilde{Y}$ in $\C$ then $X$ and $Y$ are already
isomorphic in $\D$.
\end{proof}

\begin{prop} \label{symm}
\begin{itemize}
\item[(a)]{Let $X$ and $Y$ be in $\D = D^b(H)$ for a hereditary $k$-algebra $H$. Then we have
$$\Ext^1_{\D}(Y,X) \simeq D\Ext^1_{\D}(FX,Y).$$}
\item[(b)]{Let $\widetilde{X}$ and $\widetilde{Y}$ be in $\C = \D/F$.
Then $\Ext^1_{\C}(\widetilde{X},\widetilde{Y}) \simeq  \Ext^1_{\C}(\widetilde{Y},\widetilde{X})$.}
\item[(c)]{Let $X,Y$ be indecomposable $kQ$-modules. Then
$$\Ext^1_{\C}(\widetilde{X},\widetilde{Y})\simeq \Ext^1_{kQ}(X,Y)\coprod \Ext^1_{kQ}(Y,X).$$}
\item[(d)]{If $X,Y$ are $kQ$-modules and $X$ is projective then
$$\Hom_{\C}(\widetilde{X},\widetilde{Y})\simeq \Hom_{kQ}(X,Y).$$}
\end{itemize}
\end{prop}

\begin{proof}
(a) We have $$\Ext^1_{\D}(Y,X) \simeq D \Hom_{\D}(\tau^{-1}X,Y) \simeq 
D\Ext^1_{\D}(\tau^{-1}X[1],Y) \simeq D\Ext^1_{\D}(FX,Y).$$
\noindent (b) follows directly from (a). \\
\noindent (c) We note that, by part (a), $D\Ext^1_{\D}(FX,Y)\simeq
\Ext^1_{\D}(Y,X)\simeq \Ext^1_{kQ}(Y,X)$. Suppose that $i\not=0,1$. Then
$$\Ext^1_{\D}(F^iX,Y) \simeq \Hom_{\D}(F^iX,Y[1]) \simeq 
\Hom_{\D}(F^iX,\tau FY) \simeq \Hom_{\D}(F^{i-1}\tau^{-1}X,Y).$$
The result then follows from Proposition~\ref{fewnonzero}(a), noting that
$\tau^{-1}X$ is an object in $\S$. \\
(d) If $X$ is projective, then
$$D\Hom_{\D}(F^{-1}X,Y)=D\Hom_D(\tau X[-1],Y)\simeq 
D\Hom_{\D}(\tau X,Y[1])\simeq
\Ext^1_{\D}(\tau X,Y).$$
But $\tau X\simeq I[-1]$ for some injective module $I$, so
$$D\Hom_{\D}(F^{-1}X,Y)\simeq \Ext^1_D(I[-1],Y)\simeq \Ext_D(I,Y[1])
\simeq \Ext^2(I,Y)=0.$$
The claim now follows from Proposition~\ref{fewnonzero}(a).
\end{proof}

\section{Configurations and tilting sets}

It has been shown in~\cite{mrz} that there is an interesting 
connection between cluster algebras and tilting theory for hereditary algebras.
Motivated by this, we start in this section our investigations of
tilting theory in cluster categories.

We start by recalling that (combinatorial)
$\Hom$-configurations have been investigated for the stable translation quivers
$\Z \Delta$ where $\Delta$ is a simply-laced Dynkin diagram, in connection
with the classification of the selfinjective algebras of finite representation type~\cite{rie1}.
Here a subset $\T$ of the vertices in $\Z \Delta$ is a
\emph{$\Hom$-configuration} if 
\begin{itemize}
\item[(i)]{$\Hom_{k(\Z \Delta)}(X,Y) = 0$ for all $X \neq Y$ in $\T$, and}
\item[(ii)]{for any vertex $Z$ in $\Z \Delta$ there is some $X \in \T$
such that $\Hom_{k(\Z \Delta)}(Z,X) \neq 0$.}
\end{itemize}
Of course, this can be formulated for the category $D^b(H)$ when
$\Delta$ is the underlying graph of the quiver of $H$.
$\Hom$-configurations for factors of $\Z\Delta$ are defined in the same way.

We here formulate analogous conditions using $\Ext^1$ instead
of $\Hom$, in the more general setting of the categories $\D = D^b(H)$ or $\D/G=
D^b(H)/G$ for an arbitrary finite dimensional hereditary algebra $H$. We say that a subset
$\T$ of non-isomorphic indecomposable objects in $\D$ or $\D/G$ is an
\emph{$\Ext$-configuration} if
\begin{itemize}
\item[(E1)]{$\Ext^1(X,Y) = 0$ for all $X$ and $Y$ in $\T$, and}
\item[(E2)]{for any indecomposable $Z \not \in \T$ there is some $X \in \T$
such that $\Ext^1(X,Z) \neq 0$.}
\end{itemize}
Note that in (E2) it is clearly necessary to assume that
$Z \not \in \T$.

When we have a $\Hom$-configuration $\T$ for $\Z \Delta$, with $\Delta$ Dynkin, 
it is known that $\T$ is stable under the action of $\tau^{m_{\Delta}}$.
Here $m=m_{\Delta}$ is the smallest integer such that 
in $k(\Z \Delta)$, the composition of the maps in a path of 
length greater than or equal to $m$, is zero.
Here $m_{A_n}= n$, $m_{D_n}= 2n-3, m_{E_6} = 11,m_{E_7} = 17$ 
and $m_{E_8} = 29$~\cite{blr}; in each case $m_{\Delta}=h_{\Delta}-1$,
where $h_{\Delta}$ is the Coxeter number of $\Delta$. Further, a fundamental
domain for the action of $\tau^{m_\Delta}$ has exactly $n$ objects
from $\T$, where $n$ is the number of vertices of $\Delta$, and hence
the number of non-isomorphic simple $H$-modules.

The corresponding role for $\Ext$-configurations is played by the
functor $F = \tau^{-1}[1]$ on $D^b(H)$, and this is another
reason for the importance of this functor.

\begin{prop} \label{F-closed}
Let $\T$ be an $\Ext$-configuration in $\D = D^b(H)$, and let $M$ be in $\ind \D$.
Then $M$ is in $\T$ if and only if $FM$ is in $\T$. 
\end{prop}

\begin{proof}
Assume that $M$ is in $\T$. It suffices to show that
$FM$ and $F^{-1}M$ are in $\T$. Suppose first that $F^{-1}M \not \in \T$.
Then by (E2) there is some $X$ in $\T$ such that 
$\Ext^1_{\D}(X,F^{-1}M) \neq 0$, so
$\Ext^1_{\D}(FX,M) \neq 0$. Then $\Ext^1_{\D}(M,X) \neq 0$ by
Proposition~\ref{symm},
which gives a contradiction to (E1) since $M$ and $X$ are in $\T$.
Hence we have $F^{-1}M \in \T$.

Suppose next that $FM \not \in \T$. Then by (E2) 
there is an $X$ in $\T$ such that $\Ext^1_{\D}(X,FM) \neq 0$. Since $X$ is
in $\T$ it follows that $F^{-1}X \in \T$ by the first part of the proof.
Then $\Ext^1_{\D}(F^{-1}X,M) \simeq \Ext^1_{\D}(X,FM) \neq 0$,
which contradicts (E1), since $F^{-1}X$ and $M$ are both in
$\T$.
\end{proof} 
 
There is a connection between $\Ext$-configurations in 
$\D$ and in $\D/G$, with $G$ satisfying (g1) and (g2), 
which is especially nice for $G=F$.

\begin{prop} \label{extconfiglink}
\begin{itemize}
\item[(a)]{Suppose that $\widetilde{\T}$ is an $\Ext$-configuration in the factor
category $\widetilde{D} 
= D^b(H)/G$.
Then $\T = \{X \in D^{b}(H) \mid \widetilde{X} \in \widetilde{\T} \}$ is an $\Ext$-configuration
in $\D$.}
\item[(b)]{Let $\T$ be an $\Ext$-configuration in $\D$. Then 
$\widetilde{\T} = \{ \widetilde{X} \mid X \in \T \}$ is an $\Ext$-configuration in $\C = D^b(H)/F$.}
\end{itemize}
\end{prop}

\begin{proof}
\noindent (a) Let $X$ and $Y$ be in $\T$. Then $\widetilde{X}$ and $\widetilde{Y}$ are in $\widetilde{\T}$,
so \sloppy $\Ext^1_{\widetilde{D}}(\widetilde{X}, \widetilde{Y}) = 0$. Then $\Ext^1_{\D}(X,Y) = 0$, so (E1) holds.

Let $Z \in \ind \D$, such that $Z$ is not in $\T$. Then $\widetilde{Z}$ is indecomposable in 
$\widetilde{D}$, with $\widetilde{Z} \not \in \widetilde{\T}$. So by (E2) there is
an $X \in \ind \D$ with $\widetilde{X} \in \widetilde{\T}$ such that 
$\Ext^1_{\widetilde{D}}(\widetilde{X},\widetilde{Z}) \neq 0$ 
in the factor category.
Then $\Ext^1_{\D}(G^n(X),Z) \neq 0$ for some $n$. But $G^n(X)$ lies in $\T$, since 
$\widetilde{G^n(X)} = \widetilde{X}$,
so $\T$ satisfies (E2). Hence $\T$ is an $\Ext$-configuration in $\D$.

\noindent (b) Let $X,Y$ be in $\T$, so that $\widetilde{X}, \widetilde{Y}$ are in $\widetilde{\T}$. 
Suppose for a contradiction that
$\Ext^1_{\C}(\widetilde{X}, \widetilde{Y}) \neq 0$. Then there is some integer $n$ such that 
$\Ext^1_{\D}(F^n(X),Y)) 
\neq 0$. Since $F^n(X) \in \T$ by Proposition~\ref{F-closed}, 
we have a contradiction to (E1) for $\T$.
Hence $\widetilde{\T}$ also satisfies (E1).

Now suppose that $Y \in \ind \D$ is such that $\widetilde{Y} \not \in \widetilde{\T}$. Then $Y \not \in \T$,
so there is an $X\in \T$ such that $\Ext^1_{\D}(X,Y) \neq 0$, by (E2) for $\T$. Then 
$\Ext^1_{\C}(\widetilde{X},\widetilde{Y}) \neq 0$. Since $\widetilde{X}$ is in $\widetilde{\T}$, it follows 
that $\widetilde{\T}$ 
satisfies (E2). Therefore $\widetilde{\T}$ is an $\Ext$-configuration in $\C$.
\end{proof}

The concept of $\Ext$-configurations is closely related to tilting theory for 
hereditary algebras. 
Recall that for a hereditary algebra $H$, an $H$-module $T$ is said to be a
{\em tilting module} if \\
(a) $\Ext^1_H(T,T) = 0$, that is $T$ is {\em exceptional}, and there is an exact sequence
$0 \to H \to T_0 \to T_1 \to 0$ with $T_0$ and $T_1$ in $\add T$ (see~\cite{hr}). \\
There are some useful equivalent characterisations~\cite{bo1}: \\
(b) $T$ is exceptional and has $n$
non-isomorphic indecomposable direct summands (possibly with multiplicities),
where $n$ is the number of non-isomorphic simple modules, or \\
(c) $T$ is exceptional and has a maximal number of non-isomorphic
indecomposable direct summands. \\
A tilting module is said to be {\em basic} if all of its
direct summands are non-isomorphic.

Motivated by this we say that in the categories $D^b(H)$ or $D^b(H)/G$ a set
of non-isomorphic indecomposable objects $\T$ is a \emph{tilting set} if it is
an exceptional set, that is $\Ext^1(T,T') = 0$ for all $T,T'$ in $\T$, and it
is maximal with respect to this property. 
For $D^b(H)$ there is already the concept of tilting complexes, which is
quite different, since
there the vanishing of $\Ext^i_{\D}(T,T')$ for $i \neq 0$ is required.
For the case $\C = D^b(H)/F$ we say that $T$ in $\C$ is a
{\em tilting object} if $\Ext^1_{\C}(T,T)=0$ and $T$ has a maximal number
of non-isomorphic direct summands. We note that an object in $\C$ is a basic
tilting object if and only if it is the direct sum of all objects in a
tilting set $\T$. We shall later see that all tilting sets in $\C$ are
finite, so that there will always be a corresponding basic tilting object.

We now discuss the connection between tilting sets, tilting objects and
$\Ext$-configurations.

\begin{prop} \label{induced}
Let $\widetilde{\T}$ be a set of non-isomorphic objects in $\ind\C$. Then
$\widetilde{\T}$ is a tilting set if and only if it is an
$\Ext$-configuration.
\end{prop}

\begin{proof}
Suppose that $\widetilde{\T}$ is a tilting set in $\C$. Then $\widetilde{\T}$
satisfies (E1) by definition. Let $M \in \ind \C$ such that $M \not \in \widetilde{\T}$.
If $\Ext^1_{\C}(X,M) = 0$ for all $X$ in $\widetilde{\T}$, then
$\Ext^1_{\C}(M,X) = 0$ for all $X$ in $\widetilde{\T}$ by Proposition~\ref{symm}.
Hence $\widetilde{\T} \cup \{M \}$ is exceptional, contradicting the maximality of
$\widetilde{\T}$. Hence there is some $X \in \widetilde{\T}$ such that
$\Ext^1_{\C}(X,M) \neq 0$, so that (E2) holds, so $\widetilde{\T}$ is an $\Ext$-configuration
in $\C$.

Next suppose $\widetilde{\T}$ is an $\Ext$-configuration in $\C$. Then $\widetilde{\T}$
is exceptional. By (E2), for all $M\not\in\widetilde{T}$ there is some
$X \in \widetilde{\T}$ such that $\Ext^1_{\C}(X,M) \neq 0$.
It follows that $\widetilde{\T}$ is maximal exceptional, and therefore a tilting set.
\end{proof}

Note that in $\D = D^b(H)$ there are tilting sets which are not
$\Ext$-configurations. The problem is that $\Ext^1_{\D}(\ ,\ )$ is not
symmetric.

\subsection*{Example} 
\begin{figure}[htbp]
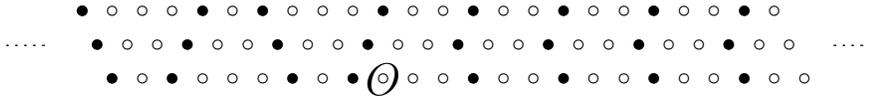

\beginpicture
\setcoordinatesystem units <0.2cm,0.45cm>             
\setplotarea x from -7 to 49, y from -1 to 3       
\linethickness=0.5pt           
\setdashes <0.2mm,1mm>
\setlinear \plot -5 1 -2 1 / %
\setlinear \plot 50 1 53 1 / %

\multiput{$\circ$}[c] at 0 2 *23 2 0 /
\multiput{$\circ$}[c] at 1 1 *23 2 0 /
\multiput{$\circ$}[c] at 2 0 *23 2 0 /

\put{$\bullet$}[c] at 0 2
\put{$\bullet$}[c] at 8 2
\put{$\bullet$}[c] at 12 2
\put{$\bullet$}[c] at 20 2
\put{$\bullet$}[c] at 26 2
\put{$\bullet$}[c] at 32 2
\put{$\bullet$}[c] at 0 2
\put{$\bullet$}[c] at 38 2
\put{$\bullet$}[c] at 44 2

\put{$\bullet$}[c] at 1 1
\put{$\bullet$}[c] at 7 1
\put{$\bullet$}[c] at 13 1
\put{$\bullet$}[c] at 19 1
\put{$\bullet$}[c] at 25 1
\put{$\bullet$}[c] at 31 1
\put{$\bullet$}[c] at 37 1
\put{$\bullet$}[c] at 43 1

\put{$\bullet$}[c] at 2 0
\put{$\bullet$}[c] at 6 0
\put{$\bullet$}[c] at 14 0
\put{$\bullet$}[c] at 18 0
\put{$\bullet$}[c] at 26 0
\put{$\bullet$}[c] at 32 0
\put{$\bullet$}[c] at 38 0
\put{$\bullet$}[c] at 44 0

\put{\huge $O$}[c] at 20 0

\endpicture
\caption{A tilting set in $\D$ which is not an $\Ext$-configuration}
\label{tiltingsetexample}
\end{figure}

\rm
Suppose that $H$ is the path algebra of a quiver of type $A_3$.
See Figure~\ref{tiltingsetexample} for the AR-quiver $\Gamma(\D)$,
indicating vertices which lie in $\T$ by filled-in circles, and those not in
$\T$ by empty circles. The arrows are omitted.
It is easy to check that $\Ext^1_{\D}(X,Y)=0$ for all $X,Y$ in $\T$, and
that $\T$ is maximal with this property, since for all
$M\not\in \T$, there is $X\in \T$ such that
$\Ext^1_{\D}(X,M)\not=0$ or $\Ext^1_{\D}(M,X)\not=0$.
In fact for all $M\not\in \T$, there is $X\in \T$ for which
$\Ext^1_{\D}(X,M)\not=0$, except for the module $N$ corresponding to the
encircled vertex. We note that $\tau N\in \T$ and
$\Ext^1_{\D}(N,\tau N)\not=0$.
So $\T$ is a tilting set in $\D$.
We note that $\T$ is not an $\Ext$-configuration, since
$\Ext^1_{\D}(X,N)=0$ for all $X\in\T$, although $N\not\in\T$.
Note also that this subset is not $F$-invariant, so could not be an
$\Ext$-configuration by Proposition~\ref{F-closed}.

We shall see in Section~\ref{relationship} that
any tilting set in $\C = D^b(H)/F$ is induced by a basic tilting module
over some hereditary algebra derived equivalent to $H$.
Hence by Proposition~\ref{induced} any $\Ext$-configuration in $\C$
is induced by such a basic tilting module. This gives another analogy with
$\Hom$-configurations, since it is known that any $\Hom$-configuration 
on $\Z \Delta$ for a Dynkin diagram $\Delta$, is induced by a basic tilting
$H$-module for a hereditary algebra $H$ whose quiver has underlying
graph $\Delta$.

Let $\Delta$ be a simply-laced Dynkin diagram, and denote by $\Pi(\Delta)$ the
preprojective algebra of type $\Delta$. Then it is known that $\Pi(\Delta)$
has finite representation type if and only if $\Delta$ is of type
$A_1$, $A_2$, $A_3$ or $A_4$ (see~\cite{dr}).
In type $A_1$, the stable module category of $\Pi(\Delta)$
has only one indecomposable (simple) object. In types $A_2$, $A_3$
and $A_4$, the stable module category of $\Pi(\Delta)$ can be seen to
coincide with the cluster category of type $A_1$, $A_3$ and $D_6$
respectively.

Let $n$ be the number of non-isomorphic simple $H$-modules, and let $t$ be the number
of non-isomorphic indecomposable $H$-modules for a hereditary algebra
of finite representation type. Then we have seen that a fundamental domain
for the action of $F$ on $D^b(H)$ has $t+n$ indecomposable objects, and we have
mentioned that there are $n$ members of an $\Ext$-configuration. For comparison,
a fundamental domain for the action of $\tau^{m_\Delta}$ is known to have $2t-n$
indecomposable objects, with $n$ members of a  $\Hom$-configuration. So we see that
in general ``more space'' is needed to have a $\Hom$-configuration. But in small
cases it may be the same, as the following example shows.

\subsection*{Example} \rm
Let $H$ be of type $A_n$. Then $t = \frac{n(n+1)}{2}$, so that we have
$\frac{n(n+1)}{2} + n = \frac{n^2 + 3n}{2}$ members of a fundamental domain for
$F$ and $n(n+1) - n = n^2$ for $\tau^{m_\Delta}$. We see that for
$n=3$, we get 9 in both cases. In this case the preprojective algebra of $H$ has
9 indecomposable nonprojective modules and induces a $\Hom$-configuration on $\Z A_3$.

\section{Relationship to tilting modules} \label{relationship}
In this section we show basic tilting modules in $\mod H$ induce tilting 
objects in $\C= D^b(H)/F$ for a hereditary algebra $H$, and
that in fact all the basic tilting objects in 
$\C$ can be obtained from basic tilting modules over
hereditary algebras derived equivalent to $H$. This allows us to deduce
additional  information on the basic tilting objects
in $\C$:  A basic exceptional object in $\C$ can be extended to a basic tilting
object, and the number of indecomposable direct summands in a 
basic tilting object is the number $n$ 
of non-isomorphic simple $H$-modules.
In particular a basic exceptional object in $\C$ with $n$ non-isomorphic
indecomposable direct summands is a basic tilting object in $\C$.

We start with the following immediate relationship between exceptional
objects in $\mod H$ and in $\C$.

\begin{lem} \label{exceptional}
Let $T$ be an $H$-module. Then $T$ is exceptional if and only if $T$ is an
exceptional object in $\C$.
\end{lem}

\begin{proof}
This follows directly from Proposition \ref{symm}(c).
\end{proof}

We use this to show the following.

\begin{prop} \label{extending}
Let $H$ be a hereditary algebra with $n$ non-isomorphic simple modules,
and let $T$ be a basic exceptional object in $\C$. Then $T$ can be extended
to a basic tilting object.
\end{prop}

\begin{proof}
We claim that any basic exceptional object $T'$ in $\C$ has at most $2n$
indecomposable summands. Let $T_1, \dots, T_r, T_{r+1}, \dots T_t$
be indecomposable objects in $\mod H \vee H[1]$ such that  
$T_1 \coprod \cdots \coprod T_r \coprod \cdots \coprod T_t$ determines
$T'$, where $T_1, \dots, T_r$ are in $\mod H$ and $T_{r+1}, \dots, T_t$ are
summands of $H[1]$. Then $T_1 \coprod \cdots \coprod T_r$ is a basic
exceptional $H$-module, so that $r \leq n$, and hence $t \leq 2n$.
In particular any basic exceptional object having $T$ as a direct summand has
at most $2n$ indecomposable direct summands, and hence $T$ can be extended to
a maximal basic exceptional object in $\C$, which is then by definition a
basic tilting object.
\end{proof}

By Lemma~\ref{exceptional}, a basic tilting $H$-module gives rise to a
basic exceptional object in $\C$ (as indecomposable $kQ$-modules are
isomorphic as modules if and only if they are isomorphic in $\C$).
We shall show that this is in fact a basic tilting object, and that any
basic tilting object in $\C$ can be obtained this way.

\begin{thm} \label{thm-sec3}
\begin{itemize}
\item[(a)]{Let $T$ be a basic tilting object in $\C = D^b(H)/F$,
where $H$ is a hereditary algebra with $n$ simple modules.
\begin{itemize}
\item[(i)]{$T$ is induced by a basic tilting module over a hereditary
algebra $H'$, derived equivalent to $H$.}
\item[(ii)]{$T$ has $n$ indecomposable direct summands.}
\end{itemize}
}
\item[(b)]{Any basic tilting module over a hereditary algebra $H$ induces
a basic tilting object for $\C = D^b(H)/F$.}
\end{itemize}
\end{thm}

\begin{proof}
(a)(i) Let $T$ be a basic tilting object in $\C = D^b(H)/F$.
Let $T_1, \dots ,T_r$ be indecomposable objects in $\mod H \vee H[1]$
inducing $T$.
If no $T_i$ is a summand of $H[1]$, then $T_1 \coprod \cdots \coprod T_r$ is
a basic exceptional $H$-module which we claim is a basic tilting module.
If not, we get a basic tilting module by adding a nonzero module as summand.
But then this will give rise to a basic exceptional object in $\C$ properly
containing $T$ as a direct summand, which is a contradiction to $T$ being a
basic tilting object in $\C$.

If no $T_i$ is projective, we have
$$\{T_1, \dots, T_r \} \subset \tau^{-1}_{\D}(\mod H)$$
and then $T_1 \coprod \cdots \coprod T_r$ is a basic tilting module over a
hereditary algebra derived equivalent to $H$
(in fact isomorphic to $H$, but with a different embedding into $D^b(H))$.
Assume now that some $T_i$ is projective.
Let first $H$ be of infinite representation type.
We assume that there are
some $T_j$ which are summands of $H[1]$ (otherwise we are done, by
the above argument).
If $T$ has no injective direct summands, then
$\tau_{\C}^{-1} T$ can be represented by a module in $\mod H$. If
$T$ has an injective direct summand (such that $\tau_{\C}^{-1} T$
has a summand in $H[1]$), we can apply $\tau_{\C}^{-1}$ again.
It is clear that there is a $t$ such that 
$\tau^{-t}_{\C} T$ can be represented by a module in $\mod H$.
Hence, $T$ is a module over 
a hereditary algebra derived equivalent to $H$, and we proceed as above.

Let now $H$ be of finite representation type, and we use the same notation as above, 
with $T_1, \dots ,T_r$ in $\mod H \vee H[1]$. We claim that for any
simple projective module $S$ not in $\add T$, there is a path to some $T_i$.
Since $T$ is a basic tilting object, we have $\Ext^1_{\C}(T,S) \neq 0$, and
hence
$\Hom_{\C}(S, \tau_{\C} T) \neq 0$. Since $\Hom_{\D}(S,F(\tau_{\D} T)) = \Hom_{\D}(S,T[1]) = 0$,
we must have $\Hom_{\D}(S,\tau_{\D} T) \neq 0$
(using Proposition~\ref{fewnonzero}(a)), and consequently we have a path of
the desired type.
Denote by $\alpha(H)$ the sum of the lengths of all paths
(where paths through the same sequence of vertices are counted only once)
from a simple projective $H$-module which is not in $\add T$, to some $T_i$.
By possibly replacing $H$ by a derived equivalent hereditary algebra, we
can assume that $\alpha(H)$ is smallest possible, when all $T_i$ are in $\mod H \vee H[1]$. 
If $\alpha(H) > 0$, there is some simple projective $H$-module $S$ not in $\add T$.
By performing an APR-tilt (see~\cite{apr}) using the basic tilting module
$M = \tau^{-1}S \coprod P$, where $H = S \coprod P$, to get
$H' = \End_H(M)^{\op}$, it is easy to see that $\alpha(H') < \alpha(H)$,
and that $H'$ satisfies the desired properties. This contradiction implies that 
$\alpha(H) = 0$, so that all simple projective $H$-modules are in $\add T$.
 
We next want to show that no $T_i$ is a summand of $\tau_{\D}^{-1} H$.
Assume to the contrary that
there is an indecomposable projective $H$-module $P$ with $\tau_{\D}^{-1}P$
in $\add T$. There is a simple projective $H$-module $S$ with
$\Hom_H(S,P) \neq 0$, and
as we have seen, it is in $\add T$. Since
$\Ext^1_{\D}(\tau^{-1}P,S) \simeq \Hom_H(S,P)$,
we have a contradiction to $T$ being exceptional. Hence no $T_i$ is a summand
of $\tau_{\D}^{-1} H$.

Choose $H'$ derived equivalent to $H$ such that
$\tau^{-2}_{\D}(\mod H \vee H[1]) = \mod H' \vee H'[1]$. Since no $T_i$ is a 
summand of $\tau_{\D}^{-1}H$, no $T_i$ is a summand of $H'[1]$ (now
regarding the $T_i$ as objects in $\mod H'\vee H'[1]$; see
Proposition~\ref{indecomposableobjects}). So $T$ is a
basic exceptional $H'$-module which has to be a basic tilting module.

(a)(ii) This is clearly a consequence of part (i).

(b) Let $T$ be a basic exceptional object in $\C$ induced by a basic tilting
$H$-module. Then
$T$ has $n$ indecomposable direct summands, and can be extended to a basic
tilting object by 
Proposition~\ref{extending}. But any basic tilting object has $n$
indecomposable direct summands, and consequently $T$ is a basic tilting
object in $\C$.
\end{proof}

We note that the basic tilting modules of $kQ$ are in bijection with the
$\Hom$-configurations of $\D$~\cite{blr}.
The above Theorem indicates a link between
tilting sets in $\C$ and basic tilting modules, which, in the light of
Propositions~\ref{extconfiglink} and~\ref{induced}, gives a link between
$\Ext$-configurations in $\D$ and basic tilting modules.
It would be interesting to find a direct link between the
$\Hom$-configurations and the $\Ext$-configurations of $\D$.

The previous investigation holds more generally 
in the setting of a hereditary abelian category $\sH$
with finite dimensional $\Hom$-spaces and $\Ext$-spaces and
with a tilting object $T$, as introduced and investigated in \cite{hrs}.
We still have Serre duality for $D^b(\sH)$ and hence almost split triangles,
(see~\cite{hrs},~\cite{rv}) and Keller's theorem on 
$\C_{\sH} = D^b(\sH)/F$ being triangulated is proved in this
generality~\cite{k}. It is also known in this setting
that a basic object $T$ in $\sH$ is a tilting object if and only if 
$\Ext^1_{\sH}(T,T) = 0$ and the number of indecomposable direct summands
of $T$ is equal to the rank of the Grothendieck-group of $\sH$. Furthermore,
any exceptional object can be extended to a tilting object, see \cite{hu1}.
Using this, the previous results carry over to this setting.

When $\sH$ is connected and not equivalent to some
$\mod H$ for a hereditary algebra $H$, it is known that
$\sH$ has no non-zero projective or injective objects, see \cite{hu2}.

In this case it is clear that
$\ind \sH$ is a fundamental domain for
$\C$ under the action of $F$. For if $X$ is in
$\ind \sH$, then $F^i X$ is in $\ind \sH[i]$, so that
no other object in the $F$-orbit of $X$ is in $\ind \sH$.
And given any $Y$ in $\ind \D$, we have
$Y[i] \in \sH$ for some $i$, and so $F^i Y=\tau^{-i}Y[i]$ is
in $\sH$ since $\sH$ is closed under positive and
negative powers of $\tau$. We then get the following.

\begin{prop} \label{corre}
Let $\sH$ be a hereditary abelian $k$-category over a field $k$,
with finite dimensional $\Hom$-spaces and $\Ext^1$-spaces.
Assume $\sH$ has no nonzero projective or injective objects, and assume
that $\sH$ has a tilting object. Then there is
a natural 1--1 correspondence between the exceptional objects
in $\sH$ and in $\C_{\sH} = D^b(\sH)/F$. The correspondence preserves
tilting objects.
\end{prop}

As has previously been done for $\mod H$ and
other hereditary categories $\sH$ with tilting objects
(see \cite{hu2}), one can associate to
$\mathcal{C}$ a tilting graph whose vertices are the basic tilting
objects, and where there is an edge between two vertices if the
corresponding tilting objects have all but one indecomposable
summands in common. It is known that for $\operatorname{mod} H$
the graph is not always connected, but this is the case for the
hereditary abelian $\Ext$-finite $k$-categories with tilting
objects derived equivalent, but not equivalent, to
$\operatorname{mod} H$ \cite{hu2}. Using this last result, we obtain
the following.

\begin{prop}
For an indecomposable hereditary
$k$-algebra $H$, the tilting graph of $\mathcal{C}=\mathcal{C}_H$
defined above is connected.
\end{prop}

\begin{proof} If $H$ is given by a Dynkin diagram, the
tilting graph for $\operatorname{mod} H$ is connected, as pointed
out in \cite{hu2}, and hence the same is true for the tilting graph of
$\mathcal{C}$.

If $H$ is of infinite representation type, it is known that there
is some indecomposable hereditary abelian
$k$-category $\sH$ with tilting objects, finite dimensional
$\Hom$-spaces and $\Ext$-spaces and no nonzero
projective or injective objects, with $D^b(\mod H)$ 
equivalent to $D^b(\sH)$ (see e.g. \cite{hu2}).
Consequently $\C_H$ is equivalent to $\C_{\sH}$.
It follows from Proposition \ref{corre} that 
the tilting graph for $\sH$ and $\C_{\sH}$ are isomorphic.
Since it is proved in \cite{hu2} that the tilting
graph of $\sH$ is connected, our result follows.
\end{proof}

\section{Connections with Cluster Algebras} \label{connections}

In this section, we assume that $H$ is the path algebra of a
simply-laced quiver of Dynkin type, with underlying graph $\Delta$, and that
$k$ is algebraically closed. We denote
by $\A=\A(\Delta)$ the corresponding cluster algebra~\cite{fz2}.
Let $\Phi$ denote the set of roots of the corresponding Lie algebra,
and let $\Phi_{\geq -1}$ denote the set of {\em almost positive} roots,
i.e.\ the positive roots together with the negatives of the simple roots.
The cluster variables of $\A$ are in 1--1 correspondence with the
elements of $\Phi_{\geq -1}$. Fomin and Zelevinsky associate a nonnegative
integer $(\alpha||\beta)$, known as the {\em compatibility degree}, to each
pair $\alpha,\beta$ of almost positive roots. This
is defined in the following way. Let $s_i$ be the Coxeter generator of
the Weyl group of $\Phi$ corresponding to $i$, and let $\sigma_i$ be the
permutation of $\Phi_{\geq -1}$ defined as follows:
$$\sigma_i(\alpha)=\left\{ \begin{array}{ll} \alpha & \alpha=-\alpha_j,\ j\not=i \\ s_i(\alpha) & \mbox{otherwise.} \end{array}\right.$$
Let $I=I^+\sqcup I^-$ be a partition of the set of vertices $I$ of
$\Delta$ into completely disconnected subsets and define:
$$\tau_{\pm}=\prod_{i\in I^{\pm}}\sigma_i.$$
Then $(\ ||\ )$ is defined by setting $(-\alpha_i||\beta)$ to be the
coefficient of $\alpha_i$ in $\beta$, and by specifying that it is
$\tau_{\pm}$-invariant.

In~\cite{mrz}, it was shown that the combinatorics of $\A$
could be obtained from the category
of {\em decorated representations} of a quiver $Q$ with underlying graph
$\Delta$.
In particular, this allowed the generalised associahedra (Stasheff polytopes)
of~\cite{cfz} to be constructed directly from the representation theory
of $Q$, and gave, for the first time, a uniform formula for the
number of basic tilting modules over $kQ$ in terms of the degrees of the
corresponding Weyl group. The compatibility degree, key to the construction
of the associahedron, was interpreted as the dimension of a certain
bifunctor from the decorated category to the category of finite-dimensional
vector spaces, in the case where the quiver was alternating.
This bifunctor can be regarded as a symmetrised version of $\Ext^1$.

In this section, we will show that such a construction can be made in a
more symmetric way, via the category $\C=D^b(H)/F$. This approach has the
advantage that the category $\C$ is independent of the orientation of the
quiver considered. We show that, when the indecomposable objects of $\C$ are
labelled appropriately with decorated representations (in a way dependent on
the orientation of the quiver) the dimension of an $\Ext^1$-group coincides
with the dimension of the symmetrised bifunctor mentioned above. Thus, when
$\C$ is labelled in a way corresponding to the alternating quiver, the
combinatorics of the corresponding cluster algebra is recovered in terms of
$\Ext^1$-groups of $\C$. In particular, we will show that the clusters are in
1--1 correspondence with the basic tilting objects in $\C$.

We first of all show that the $\Ext^1$-groups in $\C$
coincide with the symmetrised $\Ext^1$-groups used for the decorated
representations in~\cite{mrz}.
Recall that in~\cite{mrz} the quiver $Q$, with vertices $Q_0$ and arrows
$Q_1$, is replaced by a ``decorated''
quiver $\widetilde{Q}$, with an extra copy $Q_0^-=\{i_-\,:\,i\in Q_0\}$
of the vertices of $Q$ (with no arrows incident with the new copy).
A module $M$ over $k\widetilde{Q}$ can be written in the form
$M^+\coprod V$, where $M^+=\coprod_{i\in Q_0}M^+_i$ is a $kQ$-module, and
$V=\coprod_{i\in Q_0}V_i$ is a $Q_0$-graded
vector space over $k$. Its {\em signed dimension vector},
$\mathbf{sdim}(M)$ is the element of the root lattice of the Lie algebra
of type $\Delta$ given by
$$\mathbf{sdim}(M)=\sum_{i\in Q_0}\dim(M^+_i)\alpha_i-
\sum_{i\in Q_0}\dim(V_i)\alpha_i,$$
where $\alpha_1,\alpha_2,\ldots ,\alpha_n$ are the simple roots.
By Gabriel's Theorem, the indecomposable objects of $k\widetilde{Q}$-mod are
parametrised, via $\mathbf{sdim}$, by the almost positive roots,
$\Phi_{\geq -1}$, of the corresponding Lie algebra, i.e.\ the positive roots
together with the negative simple roots.
The positive roots correspond to the indecomposable $kQ$-modules, and the
the negative simple roots correspond to the simple modules associated with
the new vertices. We denote the simple module corresponding to the vertex
$i_-$ by $S_i^-$. Let $M=M^+\coprod V$ and $N=N^+\coprod W$ be two
$k\widetilde{Q}$-modules. The symmetrised $\Ext^1$-group for this pair of
modules is defined to be:
\begin{eqnarray*}
E_{kQ}(M,N) & := & \Ext^1_{kQ}(M^+,N^+)\coprod \Ext^1_{kQ}(N^+,M^+)\coprod \\
& & \Hom^{Q_0}(M^+,W) \coprod \Hom^{Q_0}(V,N^+),
\end{eqnarray*}
where $\Hom^{Q_0}$ denotes homomorphisms of $Q_0$-graded vector spaces.

We define a map $\psi_Q$ from $\ind\C$ to
the set of isomorphism classes of indecomposable $k\widetilde{Q}$-modules as
follows. Let $\widetilde{X}\in\ind\C$.
We can assume that one of the following cases holds:
\begin{enumerate}
\item $X$ is an indecomposable $kQ$-module $M^+$.
\item $X=P_i[1]$ where $P_i$ is the indecomposable projective $kQ$-module
corresponding to vertex $i\in Q_0$.
\end{enumerate}
We define $\psi_Q(\widetilde{X})$ to be $M^+$ in Case (1), and to be $S_i^-$
in Case (2).

The following is clear:

\begin{prop} \label{objectbijection}
The map $\psi_Q$ is a bijection between $\ind\C$ and the set of
isomorphism classes of indecomposable $k\widetilde{Q}$-modules
(i.e.\ indecomposable decorated representations). It follows that
$\gamma_Q:=\mathbf{sdim}\circ \psi_Q$ is a bijection between
$\ind\C$ and $\Phi_{\geq -1}$ (and thus induces a bijection between
$\ind\C$ and the set of cluster variables).
\end{prop}

For $\alpha\in \Phi_{\geq -1}$ we denote by $M_Q(\alpha)$ the element
of $\ind\C$ such that $\gamma_Q(M_Q(\alpha))=\alpha$.

\begin{prop} \label{extinterpretation}
Let $X,Y$ be objects of $\D$. Then
$$E_{kQ}(\psi_Q(\widetilde{X}),\psi_Q(\widetilde{Y}))\simeq \Ext^1_{\C}(\widetilde{X},\widetilde{Y}).$$
\end{prop}

\begin{proof}
Without loss of generality, we can assume
that $X$ and $Y$ are either indecomposable $kQ$ modules or of the form
$P_i[1]$ where $P_i$ is an indecomposable projective $kQ$-module.
We first of all consider the case where $X=M^+$ and $Y=N^+$ are both
indecomposable $kQ$-modules.
Then $E_{kQ}(\psi_Q(\widetilde{X}),\psi_Q(\widetilde{Y}))=
\Ext^1_{kQ}(M^+,N^+)\coprod \Ext^1_{kQ}(N^+,M^+)$
which is isomorphic to $\Ext^1_{\C}(\widetilde{M^+},\widetilde{N^+})$ by
Proposition~\ref{symm}.
Next, suppose that $X=P_i[1]$ and that $Y=N^+$, where $P_i$ is an
indecomposable projective and $N^+$ is an indecomposable $kQ$-module.
Then
$$E_{kQ}(\psi_Q(\widetilde{X}),\psi_Q(\widetilde{Y}))=\Hom^{Q_0}(S_i^-,N^+)$$
and has dimension given by the multiplicity of $\alpha_i$ in the positive root
corresponding to $N^+$. We also have:

\begin{eqnarray*}
\Ext^1_{\C}(\widetilde{X},\widetilde{Y}) &\simeq& \Ext^1_{\C}(\widetilde{P_i[1]},\widetilde{N^+}) \\
&\simeq& \Ext^1_{\C}(\tau^{-1}\widetilde{P_i[1]},\tau^{-1}\widetilde{N^+}) \\
&\simeq& \Ext^1_{\C}(\widetilde{P_i},\tau^{-1}\widetilde{N^+}) \\
&\simeq& \Ext^1_{\C}(\widetilde{\tau^{-1}N^+},\widetilde{P_i}) \\
&\simeq& \Hom_{\C}(\widetilde{P_i},\widetilde{N^+}) \\
&\simeq& \Hom_{kQ}(P_i,N^+),
\end{eqnarray*}
the last step by Proposition~\ref{symm}. This also has dimension equal to
the multiplicity of $\alpha_i$ in the positive root corresponding to $N^+$.

In this situation, we also have
$\Ext^1_{\C}(\widetilde{Y},\widetilde{X})\simeq \Ext^1_{\C}(\widetilde{X},\widetilde{Y})$ of the same dimension,
and $E_{kQ}(\psi_Q(\widetilde{Y}),\psi_Q(\widetilde{X}))=\Hom^{Q_0}(N^+,S_i^-)$ with
the same dimension, so the only case left to consider is when
$X=P_i[1]$ and $Y=P_j[1]$ where $P_i$ and $P_j$ are indecomposable
$kQ$-modules. In this case,
$$E_{kQ}(\psi_Q(\widetilde{X}),\psi_Q(\widetilde{Y}))=0,$$
and 
$$\Ext^1_{\C}(\widetilde{X},\widetilde{Y})=\Ext^1_{\C}(\widetilde{P_i[1]},\widetilde{P_j[1]})\simeq
\Ext^1_{\C}(\widetilde{P_i},\widetilde{P_j})=0$$ by
Proposition~\ref{symm}.
\end{proof}

This proposition shows that $\C$, which is independent of the orientation
of its defining quiver, can be regarded as a ``symmetrised'' (orientation
independent) version of the decorated categories
$k\widetilde{Q}$-mod, since $E_{kQ}$ can be modelled for all orientations
$Q$ of $\Delta$ by $\C$, via the labellings $\psi_Q$.

We therefore have:

\begin{cor} \label{Qdegree}
Let $\alpha,\beta\in\Phi_{\geq -1}$. Then we have
$$(\alpha||\beta)_Q=\dim\Ext^1_{\C}(M_Q(\alpha),M_Q(\beta)),$$
where $(\alpha||\beta)_Q$ denotes the $Q$-compatibility degree of
$\alpha$ and $\beta$ (see~\cite[Eq.(3.3)]{mrz}).
\end{cor}

We also have the following consequences.
Let $\Delta(\C)$ be the abstract simplicial complex on $\C$
with simplices given by the exceptional sets in $\C$, i.e. the subsets
of tilting sets. Thus the maximal simplicies are the tilting sets.

\begin{cor} \label{tiltingobjectcomplex}
Let $Q$ be any quiver of type $\Delta$.
Then $\Delta(\C)$ is isomorphic to the abstract simplicial complex
$\Delta_Q$ of~\cite[3.7,4.11]{mrz}.
\end{cor}

Corollary~\ref{tiltingobjectcomplex}, together
with~\cite[4.11,4.12]{mrz} show that the simplicial complex $\Delta(\Phi)$
of~\cite[p6]{fz5} can be obtained in a natural way from the
category $\C$ associated to $\Phi$.

\begin{theorem} \label{clusterbijection}
Let $Q=Q_{alt}$ be an alternating quiver of type $\Delta$.
Then the map $\alpha\mapsto M_{Q_{alt}}(\alpha)$ between $\Phi_{\geq -1}$ and
$\ind\C$ induces a bijection
between the following sets: \\
(1) The set of clusters in a cluster algebra of type $\Delta$. \\
(2) The set of basic tilting objects in $\C(\Delta)$.
\end{theorem}

\begin{proof}
The induced map is given by applying the map
$\alpha\mapsto M_{Q_{alt}}(\alpha)$
pointwise to a cluster regarded as a subset of $\Phi_{\geq -1}$.
The result follows from Corollary~\ref{Qdegree} and~\cite[4.12]{mrz}.
\end{proof}

\subsection*{Example}
\begin{figure}[htbp]
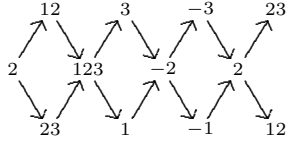


\beginpicture

\setcoordinatesystem units <0.5cm,0.4cm>             
\setplotarea x from -8 to 10, y from -9 to -3       

\scriptsize{

\put{$23$}[c] at 2 -8
\put{$1$}[c] at 4 -8
\put{$-1$}[c] at 6 -8
\put{$12$}[c] at 8 -8
\put{$2$}[c] at 1 -6
\put{$123$}[c] at 3 -6
\put{$-2$}[c] at 5 -6
\put{$2$}[c] at 7 -6
\put{$12$}[c] at 2 -4
\put{$3$}[c] at 4 -4
\put{$-3$}[c] at 6 -4
\put{$23$}[c] at 8 -4
}

\setlinear \plot  2 -7.65   2.6 -6.4  / %
\setlinear \plot  2.6 -6.4 2.7 -6.7  / %
\setlinear \plot  2.6 -6.4 2.3 -6.5  / %

\setlinear \plot  4 -7.65   4.6 -6.4  / %
\setlinear \plot  4.6 -6.4 4.7 -6.7  / %
\setlinear \plot  4.6 -6.4 4.3 -6.5  / %

\setlinear \plot  6 -7.65   6.6 -6.4  / %
\setlinear \plot  6.6 -6.4 6.7 -6.7  / %
\setlinear \plot  6.6 -6.4 6.3 -6.5  / %

\setlinear \plot  1.6 -7.65 1.0 -6.4  / %
\setlinear \plot  1.6 -7.65 1.3 -7.55  / %
\setlinear \plot  1.6 -7.65 1.7 -7.35 / %

\setlinear \plot  3.6 -7.65 3.0 -6.4  / %
\setlinear \plot  3.6 -7.65 3.3 -7.55 / %
\setlinear \plot  3.6 -7.65 3.7 -7.35  / %

\setlinear \plot  5.6 -7.65 5.0 -6.4  / %
\setlinear \plot  5.6 -7.65 5.3 -7.55  / %
\setlinear \plot  5.6 -7.65 5.7 -7.35  / %

\setlinear \plot  7.6 -7.65 7.0 -6.4  / %
\setlinear \plot  7.6 -7.65 7.3 -7.55  / %
\setlinear \plot  7.6 -7.65 7.7 -7.35  / %

\setlinear \plot  1   -5.6 1.6 -4.4  / %
\setlinear \plot  1.6 -4.4 1.7 -4.7  / %
\setlinear \plot  1.6 -4.4 1.3 -4.5  / %

\setlinear \plot  3   -5.6 3.6 -4.4  / %
\setlinear \plot  3.6 -4.4 3.7 -4.7  / %
\setlinear \plot  3.6 -4.4 3.3 -4.5  / %

\setlinear \plot  5   -5.6 5.6 -4.4  / %
\setlinear \plot  5.6 -4.4 5.7 -4.7  / %
\setlinear \plot  5.6 -4.4 5.3 -4.5  / %

\setlinear \plot  7   -5.6 7.6 -4.4  / %
\setlinear \plot  7.6 -4.4 7.7 -4.7  / %
\setlinear \plot  7.6 -4.4 7.3 -4.5  / %

\setlinear \plot  2.6 -5.65 2.0 -4.4  / %
\setlinear \plot  2.6 -5.65 2.3 -5.55  / %
\setlinear \plot  2.6 -5.65 2.7 -5.35 / %

\setlinear \plot  4.6 -5.65 4.0 -4.4  / %
\setlinear \plot  4.6 -5.65 4.3 -5.55  / %
\setlinear \plot  4.6 -5.65 4.7 -5.35 / %

\setlinear \plot  6.6 -5.65 6.0 -4.4  / %
\setlinear \plot  6.6 -5.65 6.3 -5.55  / %
\setlinear \plot  6.6 -5.65 6.7 -5.35 / %

\endpicture
\caption{The labelled AR-quiver of $\C$ in type $A_3$}
\label{clusterexample}
\end{figure}

In Figure~\ref{clusterexample}, we indicate the labelling of $\ind\C$
(via its AR-quiver) from Theorem~\ref{clusterbijection} in
type $A_3$. Objects with the same label are identified.
A positive root $\alpha_i+\alpha_{i+1}+ \cdots +\alpha_j$ is denoted
by $i,i+1,\ldots ,j$ and a negative root $-\alpha_i$ is denoted by $-i$.

\begin{figure}[htbp]
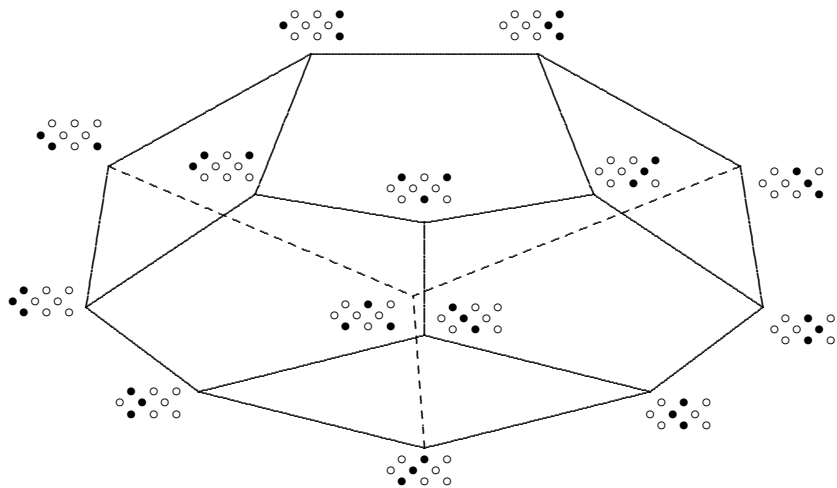


\beginpicture
\setcoordinatesystem units <0.15mm,0.075mm> point at 0 0  
\setplotarea x from -425 to 400, y from -140 to 825       
\linethickness=0.5pt           

\scriptsize{
\setlinear \plot 0 0  200 100 / %
\setlinear \plot 200 100  300 250 / %
\setlinear \plot 300 250  280 500 / %
\setlinear \plot 280 500  100 700 / %
\setlinear \plot 100 700  -100 700 / %
\setlinear \plot -100 700  -280 500 / %
\setlinear \plot -280 500  -300 250 / %
\setlinear \plot -300 250  -200 100 / %
\setlinear \plot -200 100  0 0 / %

\setlinear \plot 200 100 0 200 / %
\setlinear \plot 0 200 -200 100 / %
\setlinear \plot 0 200 0 400 / %
\setlinear \plot 0 400 150 450 / %
\setlinear \plot 150 450 300 250 / %
\setlinear \plot 150 450 100 700 / %
\setlinear \plot 0 400 -150 450 / %
\setlinear \plot -150 450 -300 250 / %
\setlinear \plot -150 450 -100 700 / %

\setdashes <1mm,1mm>
\setlinear \plot 0 0  -10 270 / %
\setlinear \plot -280 500  -10 270 / %
\setlinear \plot 280 500  -10 270 / %

\setsolid

\setcoordinatesystem point at 40 20  
\put{$\bullet$}[c] at 20 -40
\put{$\circ$}[c] at 40 -40
\put{$\circ$}[c] at 60 -40
\put{$\circ$}[c] at 10 -20
\put{$\bullet$}[c] at 30 -20
\put{$\circ$}[c] at 50 -20
\put{$\circ$}[c] at 20 0
\put{$\bullet$}[c] at 40 0
\put{$\circ$}[c] at 60 0

\setcoordinatesystem point at -190 -80  
\put{$\circ$}[c] at 20 -40
\put{$\bullet$}[c] at 40 -40
\put{$\circ$}[c] at 60 -40

\put{$\circ$}[c] at 10 -20
\put{$\bullet$}[c] at 30 -20
\put{$\circ$}[c] at 50 -20

\put{$\circ$}[c] at 20 0
\put{$\bullet$}[c] at 40 0
\put{$\circ$}[c] at 60 0

\setcoordinatesystem point at -300 -230  
\put{$\circ$}[c] at 20 -40
\put{$\bullet$}[c] at 40 -40
\put{$\circ$}[c] at 60 -40

\put{$\circ$}[c] at 10 -20
\put{$\circ$}[c] at 30 -20
\put{$\bullet$}[c] at 50 -20

\put{$\circ$}[c] at 20 0
\put{$\bullet$}[c] at 40 0
\put{$\circ$}[c] at 60 0

\setcoordinatesystem point at -290 -490  
\put{$\circ$}[c] at 20 -40
\put{$\circ$}[c] at 40 -40
\put{$\bullet$}[c] at 60 -40

\put{$\circ$}[c] at 10 -20
\put{$\circ$}[c] at 30 -20
\put{$\bullet$}[c] at 50 -20

\put{$\circ$}[c] at 20 0
\put{$\bullet$}[c] at 40 0
\put{$\circ$}[c] at 60 0

\setcoordinatesystem point at -60 -770  
\put{$\circ$}[c] at 20 -40
\put{$\circ$}[c] at 40 -40
\put{$\bullet$}[c] at 60 -40

\put{$\circ$}[c] at 10 -20
\put{$\circ$}[c] at 30 -20
\put{$\bullet$}[c] at 50 -20

\put{$\circ$}[c] at 20 0
\put{$\circ$}[c] at 40 0
\put{$\bullet$}[c] at 60 0

\setcoordinatesystem point at 135 -770  
\put{$\circ$}[c] at 20 -40
\put{$\circ$}[c] at 40 -40
\put{$\bullet$}[c] at 60 -40

\put{$\bullet$}[c] at 10 -20
\put{$\circ$}[c] at 30 -20
\put{$\circ$}[c] at 50 -20

\put{$\circ$}[c] at 20 0
\put{$\circ$}[c] at 40 0
\put{$\bullet$}[c] at 60 0

\setcoordinatesystem point at 350 -575  
\put{$\bullet$}[c] at 20 -40
\put{$\circ$}[c] at 40 -40
\put{$\bullet$}[c] at 60 -40

\put{$\bullet$}[c] at 10 -20
\put{$\circ$}[c] at 30 -20
\put{$\circ$}[c] at 50 -20

\put{$\circ$}[c] at 20 0
\put{$\circ$}[c] at 40 0
\put{$\circ$}[c] at 60 0

\setcoordinatesystem point at 375 -280  
\put{$\bullet$}[c] at 20 -40
\put{$\circ$}[c] at 40 -40
\put{$\circ$}[c] at 60 -40

\put{$\bullet$}[c] at 10 -20
\put{$\circ$}[c] at 30 -20
\put{$\circ$}[c] at 50 -20

\put{$\bullet$}[c] at 20 0
\put{$\circ$}[c] at 40 0
\put{$\circ$}[c] at 60 0

\setcoordinatesystem point at 280 -100  
\put{$\bullet$}[c] at 20 -40
\put{$\circ$}[c] at 40 -40
\put{$\circ$}[c] at 60 -40

\put{$\circ$}[c] at 10 -20
\put{$\bullet$}[c] at 30 -20
\put{$\circ$}[c] at 50 -20

\put{$\bullet$}[c] at 20 0
\put{$\circ$}[c] at 40 0
\put{$\circ$}[c] at 60 0

\setcoordinatesystem point at -5 -250  
\put{$\circ$}[c] at 20 -40
\put{$\bullet$}[c] at 40 -40
\put{$\circ$}[c] at 60 -40

\put{$\circ$}[c] at 10 -20
\put{$\bullet$}[c] at 30 -20
\put{$\circ$}[c] at 50 -20

\put{$\bullet$}[c] at 20 0
\put{$\circ$}[c] at 40 0
\put{$\circ$}[c] at 60 0

\setcoordinatesystem point at 40 -480  
\put{$\circ$}[c] at 20 -40
\put{$\bullet$}[c] at 40 -40
\put{$\circ$}[c] at 60 -40

\put{$\circ$}[c] at 10 -20
\put{$\circ$}[c] at 30 -20
\put{$\circ$}[c] at 50 -20

\put{$\bullet$}[c] at 20 0
\put{$\circ$}[c] at 40 0
\put{$\bullet$}[c] at 60 0

\setcoordinatesystem point at 90 -255  
\put{$\bullet$}[c] at 20 -40
\put{$\circ$}[c] at 40 -40
\put{$\bullet$}[c] at 60 -40

\put{$\circ$}[c] at 10 -20
\put{$\circ$}[c] at 30 -20
\put{$\circ$}[c] at 50 -20

\put{$\circ$}[c] at 20 0
\put{$\bullet$}[c] at 40 0
\put{$\circ$}[c] at 60 0

\setcoordinatesystem point at 215 -520  
\put{$\circ$}[c] at 20 -40
\put{$\circ$}[c] at 40 -40
\put{$\circ$}[c] at 60 -40

\put{$\bullet$}[c] at 10 -20
\put{$\circ$}[c] at 30 -20
\put{$\circ$}[c] at 50 -20

\put{$\bullet$}[c] at 20 0
\put{$\circ$}[c] at 40 0
\put{$\bullet$}[c] at 60 0

\setcoordinatesystem point at -145 -510  
\put{$\circ$}[c] at 20 -40
\put{$\bullet$}[c] at 40 -40
\put{$\circ$}[c] at 60 -40

\put{$\circ$}[c] at 10 -20
\put{$\circ$}[c] at 30 -20
\put{$\bullet$}[c] at 50 -20

\put{$\circ$}[c] at 20 0
\put{$\circ$}[c] at 40 0
\put{$\bullet$}[c] at 60 0

}
\endpicture
\caption{The $14$ tilting sets of $\C$ in type $A_3$}
\label{A3clusters}
\end{figure}

We recall that it is known that there is a bijection between the
clusters in type $A_n$ and the vertices of the $n$-dimensional associahedron
--- see~\cite[1.4]{cfz}.
In Figure~\ref{A3clusters}, we show the $14$ tilting sets in
type $A_3$, associated to the vertices of the $3$-dimensional associahedron
via the bijection in Theorem~\ref{clusterbijection}.
The filled-in circles indicate the elements of the tilting set; note that the
duplicated vertices of Figure~\ref{clusterexample} do not appear in these
diagrams.

\begin{prop} \label{tiltingmodules}
Given a basic $kQ$-tilting module $T$, we can write it as a direct sum
$T=\coprod_{\alpha\in S}X_{\alpha}$ where $S\subset \Phi_+$ and
$X_{\alpha}$ is the indecomposable $kQ$-module corresponding to
$\alpha\in\Phi_+$. Let $\varepsilon(T):=\coprod_{\alpha\in S} M_Q(\alpha)$.
Then $\varepsilon(T)$ is a basic tilting object of $\C$, and $\varepsilon$
defines an embedding of the set of basic tilting $kQ$-modules into the set of
basic tilting objects of $\C$.
\end{prop}

\begin{proof}
The result follows immediately from Proposition~\ref{symm}.
\end{proof}

Let $\Delta_{mod}(Q)$ denote the complex of basic exceptional $kQ$-modules.
This is an abstract simplicial complex on the set of isomorphism classes of
indecomposable $kQ$-modules, with the simplices given by the basic exceptional
$kQ$-modules. This complex was studied by C.~Riedtmann and
A.~Schofield~\cite{rs2}, and L.~Unger~\cite{u2} following a suggestion of
C.~M.~Ringel. 

\begin{cor}
Let $Q$ be any quiver of type $\Delta$.
The map $\varepsilon$ induces an embedding of $\Delta_{mod}(Q)$ into
$\Delta(\C)$.
\end{cor}

\begin{proof}
We note that $\varepsilon$ actually defines an embedding of the set of
basic exceptional $kQ$-modules into the set of exceptional sets of $\C$.
\end{proof}

For the algebra $kQ$, where $Q$ is the quiver given by $A_n$
with linear orientation, the tilting graph of the category of
finite-dimensional $kQ$-modules, as discussed in Section 3, can be 
regarded as the skeleton of a simplicial complex with simplices the faithful
basic exceptional modules. This simplicial complex is in fact the
Stasheff associahedron of dimension $n-1$, see~\cite{bk}.

\section{Complements of almost complete basic tilting objects}

Let $H$ be a finite dimensional hereditary algebra with $n$ non-isomorphic
simple modules. An $H$-module $\overline{T}$ is said to be an
\emph{almost complete basic tilting module} if
it is basic exceptional and has $n-1$ indecomposable direct
summands. Then there is automatically an indecomposable module $M$ such that
$\overline{T} \coprod M$
is a basic tilting module. Such an indecomposable module is known as a
{\em complement} to $\overline{T}$. It is known that $\overline{T}$ can be
completed to a basic tilting module in at most two different
ways~\cite{rs1,u1} and it can be done in exactly two ways if and only
if $\overline{T}$ is sincere~\cite{hu1}, that is, each simple module occurs
as a composition factor of $\overline{T}$.
We investigate the analogous concept
for the category $\C = D^b(\mod H)/ F$, and show that in this context an almost
complete basic tilting object has exactly two complements.
Hence there is a more regular behaviour in 
$\C$. Certain classes of hereditary categories exhibit a similar
behaviour~\cite{hu2}. The analogous question has been investigated for
arbitrary artin algebras~\cite{chu}.

We say that a basic exceptional object  $\overline{T}$ in $\C$ is an
\emph{almost complete basic tilting object} if there is an indecomposable object
$M$ in $\C$ such that $\overline{T} \coprod M$ is a basic tilting object.
Then we have the following main result of this section. 

\begin{thm} \label{twocomplements}
Let $H$ be a finite dimensional hereditary algebra, and $\overline{T}$ an almost 
complete basic tilting object in $\C = D^b(H)/F$. Then $\overline{T}$ 
can be completed to a basic tilting object in $\C$ in exactly two different ways.
\end{thm}

\begin{proof}
By Theorem~\ref{thm-sec3} we can assume that $\overline{T}$ is an $H$-module.
Since $\overline{T}$ is a basic exceptional $H$-module with $n-1$ non-isomorphic
direct summands, where $n$ is the number of non-isomorphic simple $H$-modules,
$\overline{T}$ is an almost complete basic tilting module over $H$.

Assume $\overline{T}$ is sincere and let 
$M_1$ and $M_2$ be the complements in $\mod H$. Since
$\overline{T} \coprod M_1$ and $\overline{T} \coprod M_2$ are
basic tilting $H$-modules, they induce basic tilting objects in $\C$ by
Theorem~\ref{thm-sec3}.
Hence $M_1$ and $M_2$ are complements to $\overline{T}$ in $\C$.
If another complement $M_3$ comes from an $H$-module,
then $\overline{T} \coprod M_3$ would be a basic exceptional $H$-module by
Lemma~\ref{exceptional} and hence a basic tilting $H$-module, which is
impossible. Let $P$ be an indecomposable projective $H$-module. Then
$\Ext^1_{\C}(P[1],\overline{T}) \simeq \Hom_{\C}(P,\overline{T})= \Hom_H(P,\overline{T}) \neq 0$ (using Proposition~\ref{symm}(d)),
since $\overline{T}$ is sincere.
Therefore $P[1]$ can not be a complement to $\overline{T}$. Hence we have exactly 
two complements when $\overline{T}$ is a sincere $H$-module.

Assume now that $\overline{T}$ is not sincere as an $H$-module,
so that there is exactly one indecomposable $H$-module which is a complement of $\overline{T}$.
It follows as above that there are no more indecomposable $H$-modules which
induce complements of $\overline{T}$ in $\C$.
Since $\overline{T}$ is not sincere, there
is an indecomposable projective $H$-module $Q$ such that $\Hom_H(Q,\overline{T}) = 0$.

Let $\G$ be the quiver of $H$, which we can assume to be a basic algebra,
and $\G'$ the subquiver obtained by removing the vertex
$e$ of $\Gamma$ corresponding to $Q$, and all arrows starting or
ending at $e$. So the corresponding path algebra $k\G'$ is isomorphic to
$H/HeH$. Then $\overline{T}$ is clearly a $k\G'$-module, and we obviously have  
$\Ext^1_{k\Gamma'}(\overline{T}, \overline{T}) = 0$ since $\Ext^1_H(\overline{T},\overline{T}) = 0$.
Since $k\Gamma'$ has $n-1$ vertices, and $\overline{T}$ has $n-1$ non-isomorphic 
indecomposable summands, $\overline{T}$ is a basic tilting module over $k\Gamma'$. Therefore $\overline{T}$
is a faithful (and hence sincere) $k\G'$-module. In particular 
$\Hom_{k\G'}(P,\overline{T}) \neq 0$
for any indecomposable projective $k\G'$-module $P$, so that $Q$
is the only indecomposable projective $H$-module with $\Hom_H(Q,\overline{T}) = 0$.

If $P[1]$, with $P$ an indecomposable projective $H$-module, is a complement
of $\overline{T}$ in $\C = D^{b}(H)/F$,
we must have $\Ext^1_{\C}(P[1], \overline{T}) = 0$, so that $\Hom_{\C}(P, \overline{T}) =0$, and
hence $\Hom_H(P,\overline{T}) = 0$. So we must have $P \simeq Q$; in particular
at most one possibility.

Conversely, if $\Hom_{\D}(Q,\overline{T})= 0$, we have $\Ext^1_{\D}(Q[1],\overline{T}) = 0$
and
$$\Ext^1_{\D}(Q[1],F\overline{T}) = \Ext^1_{\D}(Q, \tau^{-1}\overline{T}) = \Ext^1_{\D}(I[-1], \overline{T}) = 
\Ext^2_{\D}(I,\overline{T}) = 0,$$ where $I \in \mod H$. We also have 
$$\Ext^1_{\D}(Q[1],F^{-1} \overline{T}) = \Ext^1_{\D}(Q[1], \tau \overline{T}[-1]) = 0.$$
Furthermore,
$$\Ext^1_{\D}(\overline{T}, Q[1]) = \Ext^2_{\D}(\overline{T},Q) = 0,$$
and
$$\Ext^1_{\D}(\overline{T}, F^{-1}(Q)) = \Ext^1_{\D}(\overline{T}, \tau Q[-1]) = 
\Ext^1_{\D}(\overline{T}, I[-2]) = 0,$$ where $I \in \mod H$.
Hence we see that $\overline{T} \coprod Q[1]$ is a basic tilting object in $\C= D^b(H)/F$,
so that $Q[1]$ is a complement.
\end{proof}

\section{Description of complements via approximations} \label{approx}
We shall now see how, starting with a complement of an almost complete
basic tilting
object, we can construct the other one by using minimal left and right
approximations in $\C = D^b(H)/F$. This is possible since $\C$ is a
Krull-Schmidt category. We shall also use that $\C$ is in a canonical way a
triangulated category, namely the canonical functor $\D \to \C$ is a triangle
functor.

We recall the definition of minimal left and right approximations, which come
from the theory of covariantly and contravariantly finite subcategories~\cite{as}.
Suppose that $\E$ is an additive category, that $\chi$ is an additive
subcategory of $\E$, and $E$ is an object of $\E$. A map $Y\rightarrow E$
with $Y$ an object of $\chi$ is called a {\em right $\chi$-approximation} if
the induced map $\Hom_{\E}(X,Y)\rightarrow \Hom_{\E}(X,E)$ is an epimorphism
for every object $X$ of $\chi$. There is the dual notion of a
{\em left $\chi$-approximation}.
A map $f:E\rightarrow F$ in an arbitrary category $\E$ is called
{\em right minimal} if for every $g:E\rightarrow E$ such that $fg=f$,
the map $g$ is an isomorphism. Then there is the dual notion of
{\em left minimal
map}. A right (respectively, left) approximation that is also right
(respectively, left) minimal is called a {\em minimal right }(respectively,
{\em left}) {\em approximation}. 

So let as before $\overline{T}$ be an almost complete basic tilting object
in $\C$,
and let $M$ be a complement. Let $f \colon B \to M$ be a minimal 
right $\add \overline{T}$-approximation of $M$ in $\C$, and complete this map to a triangle
\begin{equation} \label{magictriangle1}
M^{\ast} \overset{g}{\rightarrow} B \overset{f}{\rightarrow} M \to M^{\ast}[1]
\end{equation}
in $\C$. 
We show in this section that $M^{\ast}$ is the second complement to 
$\overline{T}$. This can be seen as a generalisation of a result of Happel and Unger.

\begin{prop} \cite{hu1}
Let $\overline{T}$ be a sincere almost complete tilting module over
a hereditary algebra $H$. Then there are exactly two non-isomorphic 
complements $M^{\ast}$ and $M$ in $\mod H$, and an exact sequence
\begin{equation} \label{magicsequence}
0 \to M^{\ast} \to B \to M \to 0,
\end{equation}
in $\mod H$, where $B \to M$ is a minimal right $\add \overline{T}$-approximation in $\mod H$.
\end{prop}

The exact sequence (\ref{magicsequence}) gives rise
to a triangle in $\D$, and thus to a triangle in $\C$.

\begin{lem} \label{okforseq}
Assume $\overline{T}$ is an almost complete tilting object in $\C$ induced
by a sincere almost complete tilting module in $\mod H$. Then
the triangle (\ref{magictriangle1}) in $\C$ is induced by the exact
sequence (\ref{magicsequence}).
\end{lem}

\begin{proof}
We need to show that the right $\add \overline{T}$-approximation $B \to M$ in $\mod H$,
is also a right $\add \overline{T}$-approximation in $\C$.
View (\ref{magicsequence}) as a triangle in $\D$ and apply 
$\Hom_{\D}(F^{-1}\overline{T}, \ )$ to it,
to obtain an exact sequence
$$\Hom_{\D}(F^{-1}\overline{T}, B) \to \Hom_{\D}(F^{-1}\overline{T}, M)
\to \Hom_{\D}(F^{-1}\overline{T}, M^{\ast}[1]),$$
where $\Hom_{\D}(F^{-1}\overline{T}, M^{\ast}[1]) = 
\Hom_{\D}(\tau \overline{T}, M^{\ast}[2]) = 0$. Thus,
the claim follows by Proposition \ref{fewnonzero}. 
\end{proof}

To prove that $M^{\ast}$ is a second complement to $\overline{T}$ we 
use the following preliminary results.

\begin{lem}
With the above notation, we have $\Ext^1_{\C}(\overline{T},M^{\ast}) = 0 = 
\Ext^1_{\C}(M^{\ast}, \overline{T})$.
\end{lem}  

\begin{proof}
Applying $\Hom_{\C}(\overline{T},\ )$ to the triangle
$M^{\ast} \to B \to M \to M^{\ast}[1]$ we get the 
exact sequence
\begin{multline*}\Hom_{\C}(\overline{T},M^{\ast}) \to \Hom_{\C}(\overline{T},B) \overset{\Hom_{\C}(\overline{T}, f)}{\longrightarrow} 
\Hom_{\C}(\overline{T},M) \to \\
\Ext^1_{\C}(\overline{T},M^{\ast}) \to \Ext^1_{\C}(\overline{T},B).
\end{multline*}
Since $\Ext^1_{\C}(\overline{T},B) = 0$ because $B$ is in $\add \overline{T}$
and $\Ext^1_{\C}(\overline{T},\overline{T}) = 0$, and $\Hom_{\C}(\overline{T},f)$ is surjective since
$f \colon B \to M$ is a right $\add \overline{T}$-approximation, we get 
$\Ext^1_{\C}(\overline{T},M^{\ast}) = 0$. By the symmetry property 
of $\Ext^1_{\C}(\ ,\ )$, we also get $\Ext^1_{\C}(M^{\ast},\overline{T}) = 0$.
\end{proof}

\begin{lem}
The map $g \colon M^{\ast} \to B$ is a minimal left $\add \overline T$-approximation in $\C$.
\end{lem}

\begin{proof}
Apply $\Hom_{\C}(\ ,\overline{T})$ to the triangle $M^{\ast} \to B \to M \to M^{\ast}[1]$
to get the exact sequence
$$\Hom_{\C}(B,\overline{T}) \overset{\Hom_{\C}(g,\overline{T})}{\longrightarrow} \Hom_{\C}(M^{\ast},\overline{T}) \to
\Ext^1_{\C}(M,\overline{T}).$$
Since $\overline{T} \coprod M$ is a basic tilting object in $\C$, we have 
$\Ext^1_{\C}(M,\overline{T}) = 0$, and hence 
$$\Hom_{\C}(g,\overline{T}) \colon \Hom_{\C}(B,\overline{T}) \to \Hom_{\C}(M^{\ast},\overline{T})$$ 
is surjective.  So $g \colon M^{\ast} \to B$
is a left $\add \overline{T}$-approximation.

We now show that $g \colon M^{\ast} \to B$ is a left minimal map.
If it was not, then a summand $0 \to B_1$ would split off, 
where $B_1$ is a nonzero summand of $B$.
But then $B_1 \overset{\simeq}{\rightarrow} B_1$ would be
a direct summand of $f \colon B \to M$. Since $M$ is indecomposable, we would
have $M \simeq B_1$, contradicting that $B_1$ is in $\add \overline{T}$, and
that $M$ is a complement of $\overline{T}$. Our claim then follows.
\end{proof}

\begin{lem}
$M^{\ast}$ is indecomposable.
\end{lem}

\begin{proof}
Assume that $M^{\ast} = U \coprod V$ with $U$ and $V$ nonzero.
Let $f_1 \colon U \to B_1$ and $f_2 \colon U \to B_2$ be minimal left
$\add \overline{T}$-approximations, and complete the two maps
to triangles
$$U \to B_1 \to X \to U[1]$$
and
$$V \to B_2 \to Y \to V[1].$$
The direct sum of the triangles is
$$M^{\ast} \to B \to M \to M^{\ast}[1],$$
and so $M = X \coprod Y$.
Hence $X = 0$ or $Y= 0$. If $X= 0$, then $B_1 \to 0$ is a direct
summand of $f \colon B \to M$, which contradicts $f$ being right minimal.
Similarly $Y= 0$ leads to a contradiction. Hence $M^{\ast}$ is indecomposable.
\end{proof}

\begin{lem}
$M^{\ast}$ is not in $\add \overline{T}$.
\end{lem}

\begin{proof}
If $M^{\ast}$ was in $\add \overline{T}$, then $g \colon M^{\ast} \to B$ would be an 
isomorphism, and hence $M = 0$, which is a contradiction.
\end{proof}

To show that $\overline{T} \coprod M^{\ast}$ is a basic tilting object in
$\C$, it remains to show the following.

\begin{lem} \label{Mastexceptional}
$\Ext^1_{\C}(M^{\ast},M^{\ast}) = 0$.
\end{lem}

\begin{proof}
Consider again the triangle
$$M^{\ast} \overset{g}{\rightarrow} B \overset{f}{\rightarrow} M
\to M^{\ast}[1].$$
Apply $\Hom_{\C}(\ ,M)$ to get the exact sequence
$$\Hom_{\C}(B,M) \overset{\Hom_{\C}(g,M)}{\longrightarrow}
\Hom_{\C}(M^{\ast},M) \to \Ext^1_{\C}(M,M).$$
Since $\Ext^1_{\C}(M,M) = 0$, the map $\Hom_{\C}(g,M)$ is surjective.
Hence any map $h \colon M^{\ast} \to M$ factors through $g \colon M^{\ast} \to B$. 
Now apply $\Hom_{\C}(M^{\ast},\ )$ to the triangle to get the exact sequence 
$$\Hom_{\C}(M^{\ast},B) \overset{\Hom_{\C}(M^{\ast}, f)}{\longrightarrow} \Hom_{\C}(M^{\ast},M) \to 
\Ext^1_{\C}(M^{\ast},M^{\ast}) \to \Ext^1_{\C}(M^{\ast},B),$$
where the last term is zero. To show that
$\Ext^1_{\C}(M^{\ast},M^{\ast}) = 0$ it is therefore enough to show that
$\Hom_{\C}(M^{\ast},f) \colon \Hom_{\C}(M^{\ast},B) \to \Hom_{\C}(M^{\ast},M)$
is surjective, that is, any map $h \colon M^{\ast} \to M$ 
factors through $f \colon B \to M$.
Then consider the commutative diagram in Figure~\ref{commdiag1},
\begin{figure}[htbp]
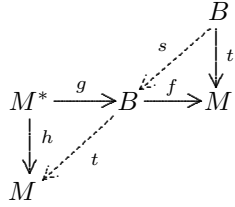

\beginpicture
\setcoordinatesystem units <0.6cm,0.6cm>             
\setplotarea x from -8 to 5, y from 0 to 4       

\put{$M$}[c] at -0.15 0
\put{$M^{\ast}$}[c] at 0 2
\put{$B$}[c] at 2.2 2
\put{$M$}[c] at 4.2 2
\put{$B$}[c] at 4.2 4
\scriptsize{
\put{$h$}[c] at 0.25 1.2
\put{$g$}[c] at 1 2.3
\put{$t$}[c] at 1.3 0.7
\put{$s$}[c] at 2.8 3.2
\put{$t$}[c] at 4.3 3
\put{$f$}[c] at 3 2.3
}
\arrow <1.8mm> [.5,1] from -0.3 1.6 to -0.3 0.4
\arrow <1.8mm> [.5,1] from 0.25 2 to 1.6 2
\arrow <1.8mm> [.5,1] from 2.25 2 to 3.5 2
\arrow <1.8mm> [.5,1] from 3.83 3.6 to 3.83 2.3

\setdashes <0.5mm,0.5mm>
\arrow <1.8mm> [.5,1] from 1.55 1.65 to 0 0.25
\arrow <1.8mm> [.5,1] from 3.7 3.6 to 2.15 2.25
\endpicture
\caption{Commutative diagram for the proof of Lemma~\ref{Mastexceptional}}
\label{commdiag1}
\end{figure}
\noindent
where $t$ is obtained from the first lifting, and we get $s \colon B \to B$ by
using that $f \colon B \to M$ is a right $\add \overline{T}$-approximation. So $h = t g= fsg$,
and hence $h \colon M^{\ast} \to M$ factors through $f \colon B \to M$, as desired.
This finishes the proof of the lemma.
\end{proof}

We now put the lemmas together to get the following.

\begin{thm} \label{magictriangle}
If $M$ is a complement of the almost complete basic tilting object
$\overline{T}$ in $\C$, then $M^{\ast}$ is another complement, obtained by
extending  the minimal right $\add \overline{T}$-approximation to a triangle.
\end{thm}
 
\begin{proof}
We only need to remark that $M \not \simeq M^{\ast}$. This follows
since $\Ext^1_{\C}(M, M^{\ast}) \neq 0$ and $\Ext^1_{\C}(M, M) = 0$.
\end{proof}

It is clear that we can also get dual constructions. That is, start with
a complement $M$, and consider the triangle 
\begin{equation}\label{magictriangle2}
M \overset{u}{\rightarrow} B' \overset{v}{\rightarrow} M^{\ast \ast} \to M[1],
\end{equation}
where $u \colon M \to B'$ is a minimal left $\add \overline{T}$-approximation.
In a dual way we get that $v \colon B' \to M^{\ast \ast}$ is a minimal
right $\add \overline{T}$-approximation, and that $M^{\ast \ast}$ is a complement of
$\overline{T}$ with $M \not \simeq M^{\ast \ast}$. We then have the following consequence of
the previous results. 
 
\begin{prop}
Let $M$ be a complement of the almost complete basic tilting object
$\overline{T}$ in $\C$. Then 
$M^{\ast} \simeq M^{\ast \ast}$ is the unique other complement, where $M^{\ast}$
is the fibre of the minimal right $\add \overline{T}$-approximation of $M$ in $\C$, and 
$M^{\ast \ast}$ is the cofibre of the minimal left $\add \overline{T}$-approximation of $M$ in $\C$.
\end{prop}

For an indecomposable exceptional module $M$, it is well known that
the endomorphism ring $\End_H(M)$ is a division ring.
However, the endomorphism ring $\End_{\C}(M)$ 
need not be a division ring, which we later in this section observe in an example.
However, if $H$ is of finite representation type, or 
more generally, if $M$ is (induced by) a preprojective or preinjective module,
then $\End_{\C}(M)$ is a division ring. This is a special case of the following.

\begin{lem}\label{preproj}
Let $M$ be an indecomposable $H$-module with $\Hom_{\D}(M, \tau^2 M) = 0$. 
Then $\End_{\C}(M)$ is a division ring.
\end{lem}

\begin{proof}
By Proposition \ref{fewnonzero}, $\End_{\C}(M) = \Hom_{\D}(M,M) \oplus \Hom_{\D}(M,FM)$.
Using the AR-formula and the assumption on $M$, we obtain $$\Hom_{\D}(M,FM)
= \Hom_{\D}(M, \tau^{-1}M[1]) \simeq D \Hom_{\D}(M, \tau^2 M) = 0$$ and
the claim follows.   
\end{proof}

We will need to consider the factor rings
$D_M = \End_{\C}(M) / \rad_{\C}(M,M)$
and $D_{M^{\ast}} = \End_{\C}(M^{\ast}) / \rad_{\C}(M^{\ast},M^{\ast})$,
which turn out to be isomorphic.

\begin{lem} \label{samefactors}
There is a natural ring-isomorphism $D_M \to D_{M^{\ast}}$.
\end{lem}

\begin{proof}
Consider the triangle
$$M^{\ast} \overset{g}{\rightarrow} B \overset{f}{\rightarrow} M \to M^{\ast}[1].$$
Let $\alpha$ be an element in $\End_{\C}(M)$. Then
there is a commutative diagram as in Figure \ref{commdiag-extra1}
\begin{figure}[htbp]
\beginpicture
\setcoordinatesystem units <0.8cm,0.6cm>             
\setplotarea x from -5 to 7, y from 0 to 3       

\put{$M^{\ast}$}[c] at 0 0
\put{$B$}[c] at 2 0
\put{$M$}[c] at 4 0
\put{$M[1]$}[c] at 6 0
\put{$M^{\ast}$}[c] at 0 2
\put{$B$}[c] at 2 2
\put{$M$} at 4 2
\put{$M[1]$} at 6 2
\scriptsize{
\put{$g$} at 0.8 0.3
\put{$f$} at 2.8 0.3
\put{$g$} at 0.8 2.3
\put{$f$} at 2.8 2.3
\put{$\beta$}[c] at 0.1 1
\put{$\gamma$}[c] at 2.2 1
\put{$\alpha$}[c] at 4.2 1
\put{$\beta[1]$}[c] at 6.4 1
}
\arrow <1.8mm> [.5,1] from 0.35 0.03 to 1.3 0.03
\arrow <1.8mm> [.5,1] from 2.25 0.03 to 3.3 0.03
\arrow <1.8mm> [.5,1] from 4.25 0.03 to 5.15 0.03
\arrow <1.8mm> [.5,1] from 0.35 2.03 to 1.3 2.03
\arrow <1.8mm> [.5,1] from 2.25 2.03 to 3.3 2.03
\arrow <1.8mm> [.5,1] from 4.25 2.03 to 5.15 2.03

\arrow <1.8mm> [.5,1] from 3.8 1.6 to 3.8 0.4
\setdashes <0.5mm,0.5mm>
\arrow <1.8mm> [.5,1] from -0.3 1.6 to -0.3 0.4
\arrow <1.8mm> [.5,1] from 1.8 1.6 to 1.8 0.4
\arrow <1.8mm> [.5,1] from 5.8 1.6 to 5.8 0.4
\endpicture
\caption{Commutative diagram for the proof of Lemma~\ref{samefactors}}
\label{commdiag-extra1}
\end{figure}
where $\gamma$ exists since $B \to M$ is a right
$\add \overline{T}$-approximation.
We claim that the map $\alpha \mapsto \beta$ gives a well-defined
ring-homomorphism
$$\End_{\C}(M) \to \End_{\C}(M^{\ast})/ \rad_{\C}(M^{\ast},M^{\ast}).$$
First note that if $g= 0$, then $B= 0$, and the map $M \to M^{\ast}[1]$ is an
isomorphism. Thus, in this case, the map $\alpha \mapsto \beta$ is well-defined.
Assume then that $g$ is non-zero.
Let $\alpha\in\End_{\C}(M)$ and fix a map $\gamma$, such that
$\alpha f = f \gamma $. Then there is some map 
from $M^{\ast}$ to $M^{\ast}$ completing the diagram.
Assume there are two such maps $\beta_1$ and $\beta_2$.
Then $g (\beta_1 - \beta_2)=0$, so $\beta_1 - \beta_2$ is not an isomorphism, and
thus each choice of $\gamma$ gives a well defined element $\overline{\beta} \in 
\End_{\C}(M^{\ast})/ \rad_{\C}(M^{\ast},M^{\ast})$. 
Let $\gamma_1$ and
$\gamma_2$ be  maps from $B$ to $B$ making the diagram commute,
and choose corresponding maps  $\beta_1, \beta_2 \in 
\End_{\C}(M^{\ast})$.
It then suffices to show that $\overline{\beta_1}-\overline{\beta_2}$ is zero, in other words
that $\beta_1 - \beta_2$ is a nonisomorphism.
We have $\alpha f = f \gamma_1 = f \gamma_2$.
Since $f(\gamma_1 - \gamma_2) = 0$, there is a map $w\colon B\rightarrow M^{\ast}$
such that $gw = \gamma_1 - \gamma_2$ and thus
$$gwg = (\gamma_1 - \gamma_2) g = g (\beta_1 - \beta_2).$$
Since $M^{\ast}$ is not a summand of $B$, $wg$ is not an isomorphism.
If $\beta_1-\beta_2$ was an isomorphism, then $wg - (\beta_1 - \beta_2)$
would also be an isomorphism.
But $$g(wg - (\beta_1 - \beta_2)) = 0,$$ so then $g=0$,
a contradiction, so the claim is proved.

Using that also $M^{\ast} \to B$ is a left $\add \overline{T}$-approximation, we obtain that
$\alpha \mapsto \overline{\beta}$ is an epimorphism. Assume now $\alpha$ is not an isomorphism, then
there is an integer $N$ such that $\alpha^N =0$ by Proposition \ref{krullschmidt}.
Thus, there is a commutative diagram as in Figure \ref{commdiag-extra2}
\begin{figure}[htbp]
\beginpicture
\setcoordinatesystem units <0.8cm,0.6cm>             
\setplotarea x from -6.2 to 5, y from 0 to 3       

\put{$B$}[c] at 0 0
\put{$M$}[c] at 2 0
\put{$M^{\ast}[1]$}[c] at 4 0
\put{$B$}[c] at 0 2
\put{$M$}[c] at 2 2
\put{$M^{\ast}[1]$}[c] at 4 2
\scriptsize{
\put{$\gamma^N$}[c] at 0.4 1
\put{0}[c] at 2.2 1
\put{$\beta^N$}[c] at 4.4 1
}
\arrow <1.8mm> [.5,1] from 0.35 0.03 to 1.3 0.03
\arrow <1.8mm> [.5,1] from 2.25 0.03 to 3.1 0.03
\arrow <1.8mm> [.5,1] from 0.35 2.03 to 1.3 2.03
\arrow <1.8mm> [.5,1] from 2.25 2.03 to 3.1 2.03

\arrow <1.8mm> [.5,1] from 3.8 1.6 to 3.8 0.4
\arrow <1.8mm> [.5,1] from -0.2 1.6 to -0.2 0.4
\arrow <1.8mm> [.5,1] from 1.8 1.6 to 1.8 0.4
\endpicture
\caption{Commutative diagram for the proof of Lemma~\ref{samefactors}}
\label{commdiag-extra2}
\end{figure}
which shows that $\beta^N$ is not an isomorphism, and thus $\beta$ is in $\rad_{\C}(M^{\ast}, M^{\ast})$.
It follows from the minimality of $B \to M$ that if $\alpha$ is an isomorphism,
then $\gamma$ and hence $\beta$ are isomorphisms.
\end{proof}

We want to show that in $\C$ all non-isomorphisms $M \to M$, actually
factor through $B \to M$. The following is useful for this.

\begin{lem} \label{symmetric}
All maps in $\rad_{\C}(M,M)$ factor through $B \to M$ if and only if
all maps in $\rad_{\C}(M^{\ast},M^{\ast})$ factor through $M^{\ast} \to B$.
\end{lem}

\begin{proof}
Apply $\Hom_{\C}(M,\ )$ to the triangle $$M^{\ast} \to B \to M \to M^{\ast}[1]$$ 
to obtain the exact sequence
$$\Hom_{\C}(M,M^{\ast}) \to \Hom_{\C}(M,B) \to \Hom_{\C}(M,M) \to 
\Hom_{\C}(M,M^{\ast}[1]) \to 0. $$  
Assume any $f \in \rad_{\C}(M,M)$ factors through $B \to M$. This means
that $$\Hom_{\C}(M,M^{\ast}[1]) \simeq \Hom_{\C}(M,M)/ \rad_{\C}(M,M).$$
Applying $\Hom_{\C}(\ , M^{\ast})$ to the same triangle gives the
exact sequence
$$\Hom_{\C}(M,M^{\ast}) \to \Hom_{\C}(B, M^{\ast}) \overset{u}{\rightarrow} 
\Hom_{\C}(M^{\ast},M^{\ast}) \to 
\Hom_{\C}(M,M^{\ast}[1]) \to 0.$$ 
which means that $\Hom_{\C}(M,M^{\ast}[1]) \simeq \Hom_{\C}(M^{\ast},M^{\ast})/ I$ 
where $I$ is the image of the map $u$.
It follows from Lemma \ref{samefactors} that $I$ 
is the radical $\rad_{\C}(M^{\ast},M^{\ast})$.
The other implication can be shown similarly.
\end{proof}

We can now prove the promised result about lifting non-isomorphisms in $\C$.

\begin{lem} \label{lifts}
With the previous notation and assumptions, 
any non-isomorphism $M \to M$ in $\C$ factors through $B \to M$.
\end{lem}

\begin{proof}
We can assume that $\overline{T} \oplus M$ is induced by an $H$-module.
We first assume that $\overline{T}$ is sincere, so $M^{\ast}$ is also
induced by a module. By Lemma \ref{okforseq}, the triangle 
(\ref{magictriangle1})
is induced by an exact sequence of modules. We view this exact sequence
as a triangle in $\D$ and apply $\Hom_{D}(F^{-1}M,\ )$ to it. 
Since $\Hom_{D}(F^{-1}M, M^{\ast}[1]) = \Hom_{D}(\tau M, M^{\ast}[2]) = 0$, it follows
that any non-isomorphism $M \to M$ in $\C$ factors through $B \to M$,
using Proposition \ref{fewnonzero}.

Now assume that $\overline{T}$ is not sincere. Then $M^{\ast}$ is induced
by an object $P[1]$ in $\D$, where $P$ is an indecomposable projective $H$-module,
and thus $\rad_{\C}(M^{\ast},M^{\ast}) = 0$, using Lemma \ref{preproj}. 
Applying Lemma \ref{symmetric}, it follows
trivially that also in this case, all non-isomorphisms $M \to M$ in $\C$
factor through $B \to M$.

\end{proof}

We can now conclude with the following property of $\Ext^1_{\C}(M,M^{\ast})$. 

\begin{prop} \label{dimensionone}
Let $M$ and $M^{\ast}$ be the complements of an almost complete tilting
object in the cluster category $\C$. Then
$\Ext^1_{\C}(M,M^{\ast})$ has dimension one over
each of the division rings $D_M = \End_{\C}(M)/ \rad_{\C}(M,M)$ and $D_{M^{\ast}}$.
\end{prop}

\begin{proof}
Apply $\Hom_{\C}(M,\ )$ to the triangle
$$M^{\ast} \to B \to M \to M^{\ast}[1],$$
to get the exact sequence
$$\Hom_{\C}(M,M^{\ast}) \to \Hom_{\C}(M,B) \overset{u}{\rightarrow} 
\Hom_{\C}(M,M) \to \Ext^1_{\C}(M,M^{\ast}) \to \Ext^1_{\C}(M,B)$$  
where $\Ext^1_{\C}(M,B) = 0$. 
Isomorphisms $M \to M$ do not lift to $B$, since $M$ is not a summand of $B$.
Thus, it follows from Lemma \ref{lifts} that
$\Ext^1_{\C}(M,M^{\ast}) \simeq \End_{\C}(M)/ \rad_{\C}(M,M)$.
It follows similarly that $\Ext^1_{\C}(M,M^{\ast})$ is one-dimensional over
$\End_{\C}(M^{\ast})$.
\end{proof}

Note that for the triangle
$$M^{\ast} \to B \to M \to M^{\ast}[1]$$
it may happen
that $B$ is zero, even though $\Ext^1_{\C}(M,M^{\ast}) \neq 0$.
This of course means that $M \simeq M^{\ast}[1] = \tau M^{\ast}$, so in
this case the second triangle $\tau M^{\ast} = M \to B' \to M^{\ast}$ is
almost split.

We also notice that Lemma \ref{lifts} has the following 
interpretation.

\begin{cor}
Let $T$ be a tilting object in a cluster category $\C$ for
a hereditary algebra algebra over an algebraically closed field.
Then the quiver of $\End_{\C}(T)^{\op}$ has no loops.
\end{cor}

\subsection*{Example}
\begin{figure}[htbp]
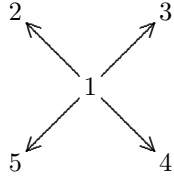

\beginpicture
\setcoordinatesystem units <1cm,1cm>             
\setplotarea x from -6.35 to 10, y from -0.5 to 2.5       

\put{$1$}[c] at 0 0
\put{$2$}[c] at -1 1
\put{$3$}[c] at 1 1
\put{$4$}[c] at 1 -1
\put{$5$}[c] at -1 -1

\arrow <1.8mm> [.5,1] from 0.15 0.15 to 0.85 0.85
\arrow <1.8mm> [.5,1] from 0.15 -0.15 to 0.85 -0.85
\arrow <1.8mm> [.5,1] from -0.15 0.15 to -0.85 0.85
\arrow <1.8mm> [.5,1] from -0.15 -0.15 to -0.85 -0.85

\endpicture
\caption{A quiver of type $\widetilde{D_4}$}
\label{d4tildequiver}
\end{figure}

The following example illustrates Theorem~\ref{magictriangle} in the tame
hereditary case. We consider the quiver $\widetilde{D}_4$ with the
orientation as in Figure~\ref{d4tildequiver}.
Let $\Lambda$ be the path algebra of the above quiver over some field $k$.
Let $P_i$ be the indecomposable projective corresponding to vertex $i$.
Then it easy to see that 
$$\overline{T} = P_4 \coprod P_5 \coprod \tau^{-1}P_2 \coprod \tau^{-1}P_3$$
is an almost complete basic tilting module, and thus an almost complete basic
tilting object in the corresponding category $\C$. 

It is clear that $\overline{T}$ is sincere, and one complement is
easy to find, namely $P_1$.
We use the above approach to find the other complement.

Let $R$ be the cokernel of the embedding 
$P_1 \to \tau^{-1} P_2 \coprod \tau^{-1} P_3$.
Then $R$ is a regular exceptional module with 
composition factors $S_1, S_4, S_5$. 
In the AR-quiver it is at the mouth of a tube of rank two, 
so $\tau^2 R \simeq R$. Thus, $R$ is an example of an 
indecomposable exceptional 
object, with $\End_{\C}(R)$ \emph{not} a division ring.

In $\mod \Lambda$ there are exact sequences
$$0 \to P_1 \to \tau^{-1} P_2 \coprod \tau^{-1} P_3 \to R \to 0$$
and
$$0 \to P_4 \coprod P_5 \to P_1 \to \tau R \to 0.$$
Thus, in $\D$ there are triangles
$$P_1 \to \tau^{-1} P_2 \coprod \tau^{-1} P_3 \to R \to P_1[1]$$
and
$$F^{-1}R \to P_4 \coprod P_5 \to P_1 \to F^{-1}R[1].$$
The images of these triangles in $\C$ are exactly 
the triangles described in Theorem~\ref{magictriangle}.

Thus, we obtain that the other complement of $\overline{T}$ is $R$, and 
$B$ and $B'$ are given by 
$P_4 \coprod P_5$ and $\tau^{-1} P_2 \coprod \tau^{-1} P_3$ respectively.

If we let $T = \overline{T} \coprod P_1$ and $T' = \overline{T} \coprod R$, then the 
endomorphism ring $\End_{\C}(T)^{\op}$ is the path algebra of the
quiver in Figure~\ref{quiver1},
while $\End_{\C}(T')^{\op}$ is the path algebra of the quiver in
Figure~\ref{quiver2}, with relations 
$ac-bf, ec-df, ga-he, gb-hd, cg, ch, fg, fh$. 
\begin{figure}[htbp]
\beginpicture
\setcoordinatesystem units <1cm,1cm>             
\setplotarea x from -5.4 to 10, y from -0.5 to 2.5       

\put{$\bullet$}[c] at 0 0
\put{$\bullet$}[c] at 2 0
\put{$\bullet$}[c] at 0 2
\put{$\bullet$}[c] at 2 2
\put{$\bullet$}[c] at 1 1

\arrow <1.8mm> [.5,1] from 0.15 0.15 to 0.85 0.85
\arrow <1.8mm> [.5,1] from 0.15 1.85 to 0.85 1.15
\arrow <1.8mm> [.5,1] from 1.15 1.15 to 1.85 1.85
\arrow <1.8mm> [.5,1] from 1.15 0.85 to 1.85 0.15

\endpicture
\caption{The quiver of $\End_{\C}(T)^{\op}$}
\label{quiver1}
\end{figure}
\begin{figure}[htbp]
\beginpicture
\setcoordinatesystem units <1cm,1cm>             
\setplotarea x from -4.3 to 10, y from -1 to 3.3       

\scriptsize{
\put{$a$}[l] at 1 2.15
\put{$b$}[l] at 0.7 1.7
\put{$c$}[l] at 3.2 1.7
\put{$d$}[l] at 1 0.15
\put{$e$}[l] at 0.4 0.7
\put{$f$}[l] at 3.1 0.7
\put{$g$}[l] at 2.2 2.75
\put{$h$}[l] at 2.2 -0.75
}

\put{$\bullet$}[c] at 0 0
\put{$\bullet$}[c] at 2 0
\put{$\bullet$}[c] at 0 2
\put{$\bullet$}[c] at 2 2
\put{$\bullet$}[c] at 4 1

\arrow <1.8mm> [.5,1] from 0.2 0 to 1.8 0
\arrow <1.8mm> [.5,1] from 0.2 2 to 1.8 2
\arrow <1.8mm> [.5,1] from 0.2 0.2 to 1.8 1.8
\arrow <1.8mm> [.5,1] from 0.2 1.8 to 1.8 0.2
\arrow <1.8mm> [.5,1] from 2.2 0 to 3.8 0.9
\arrow <1.8mm> [.5,1] from 2.2 2 to 3.8 1.1

\setlinear \plot 0.2 2.1 0.4 2.1 /
\setlinear \plot 0.2 2.1 0.25 2.28 /
\setlinear \plot 0.2 -0.1 0.4 -0.1 /
\setlinear \plot 0.2 -0.1 0.25 -0.28 /
\setquadratic
\plot 3.9 1.2 2 2.6  0.2 2.1 /
\plot 3.9 0.8 2 -0.6 0.2 -0.1 /
\endpicture
\caption{The quiver of $\End_{\C}(T')^{\op}$}
\label{quiver2}
\end{figure}
\section{Description of exchange pairs}
As usual let $H$ be a hereditary finite dimensional algebra, and
$\C$ the factor category $\D^b(H)/F$, with $F= \tau^{-1}[1]$.
We say that two non-isomorphic indecomposable objects in $\C$ form an 
\emph{exchange pair} if they are complements of the same almost complete
basic tilting object.
In this language, we have seen that if $M$ and $M^{\ast}$ form an exchange
pair, then $\Ext^1_{\C}(M,M^{\ast}) \simeq \Ext^1_{\C}(M^{\ast},M)$ is one-dimensional
over $D_M = \End_{\C}(M)/ \rad_{\C}(M,M)$ and over $D_{M^{\ast}}$. We now 
want to show that also the converse holds.

Assume that $M,M^{\ast}$ are exceptional
and that $\Ext^1_{\C}(M,M^{\ast}) = \Ext^1_{\C}(M^{\ast},M)$ is one-dimensional
over $D_M$ and over $D_{M^{\ast}}$. We can therefore choose
non-split triangles:
\begin{equation*}
M^{\ast} \to B \to M \to M^{\ast}[1]
\tag{1}\end{equation*}
and
\begin{equation*}
M \to B' \to M^{\ast} \to M[1]
\tag{2}\end{equation*}
in $\C$, where we use the same notation as before.
We want to find an almost complete basic tilting object $\overline{T}$ having
$M$ and $M^{\ast}$ as complements.
We start building up $\overline{T}$ by showing that $B \coprod B' \coprod M$
and $B \coprod B' \coprod M^{\ast}$ are exceptional objects in $\C$.

\begin{lem}
In the above notation we have:
$$\Ext^1_{\C}(B \coprod B'\coprod M,B \coprod B'\coprod M)=0$$
and
$$\Ext^1_{\C}(B \coprod B'\coprod M^{\ast},B \coprod B'\coprod M^{\ast})=0.$$
\end{lem}

\begin{proof}
Apply $\Hom_{\C}(M,\ )$ to (1) to get the exact sequence
\begin{multline*}
\Hom_{\C}(M,M^{\ast}) \to \Hom_{\C}(M,B) \to \Hom_{\C}(M,M) \overset{\alpha}{\rightarrow} \\
\Ext^1_{\C}(M,M^{\ast}) \to \Ext^1_{\C}(M,B) \to \Ext^1_{\C}(M,M).
\end{multline*}

Since $\alpha \neq 0$ and $\dim_{D_M}\Ext^1_{\C}(M,M^{\ast}) = 1$,
while
$\Ext^1_{\C}(M,M)=0$ by assumption, it follows that
$\Ext^1_{\C}(M,B)=0$. Analogously, we get $\Ext^1_{\C}(M^{\ast},B')=0$. 

Apply $\Hom_{\C}(\ ,M^{\ast})$ to (1) to get the exact sequence
\begin{multline*}
\Hom_{\C}(M,M^{\ast}) \to \Hom_{\C}(B, M^{\ast}) \to \Hom_{\C}(M^{\ast},M^{\ast}) 
\overset{\beta}{\rightarrow} \\
\Ext^1_{\C}(M,M^{\ast}) \to \Ext^1_{\C}(B, M^{\ast}) \to \Ext^1_{\C}(M^{\ast},M^{\ast}).
\end{multline*}
Since $\beta \neq 0$ and $\dim_{D_{M^{\ast}}}\Ext^1_{\C}(M,M^{\ast}) = 1$,
while
$\Ext^1_{\C}(M^{\ast},M^{\ast})=0$ by assumption, we get
$\Ext^1_{\C}(B, M^{\ast})=0$. Analogously, we get from (2) that $\Ext^1_{\C}(B',M)=0$. 

Apply $\Hom_{\C}(B \coprod B', \ )$ to (1) to get the exact sequence
$$
\Ext^1_{\C}(B \coprod B', M^{\ast}) \to \Ext^1_{\C}(B \coprod B', B) \to
\Ext^1_{\C}(B \coprod B', M),
$$
and hence $\Ext^1_{\C}(B \coprod B', B)=0$.
Apply $\Hom_{\C}(B \coprod B', \ )$ to (2) to get the exact sequence
$$\Ext^1_{\C}(B \coprod B', M) \to \Ext^1_{\C}(B \coprod B', B') \to
\Ext^1_{\C}(B \coprod B', M^{\ast}),$$ and hence $\Ext^1_{\C}(B \coprod B', B')=0$.
This finishes the proof of the lemma.    
\end{proof}

We remark that this implies that $M$ and $M^{\ast}$ cannot be direct summands
of $B\coprod B'$.
We have that $B \coprod B'$ is an exceptional object in $\C$, 
and hence can be extended to a tilting object by
Lemma~\ref{extending}.
So let $T'$ be a complement in $\C$, that is $T= B \coprod B' \coprod T'$
is a tilting object in $\C$.
We want to show that either $M$ or $M^{\ast}$ is a direct summand 
of $T$ and if we remove all copies of this summand, we get a new tilting
object by adding the other one. 

The proof of this is based upon the following crucial result. Here, for
$X$ an object of $\C$, $Supp_{\C}(\ ,X)$ denotes the objects in $\C$ which
have a non-zero homomorphism to $X$.

\begin{lem} \label{supp}
With the above notation, we have
$$\Supp_{\C}(\ ,\tau M) \subset \{M^{\ast}\} \cup \Supp_{\C}(\ ,\tau B) \cup \Supp_{\C}(\ ,\tau B').$$ 
\end{lem}

\begin{proof}
Consider the triangles $$M^{\ast} \overset{g}{\rightarrow} B \overset{f}{\rightarrow} 
M \to M^{\ast}[1]$$ and
$$M \to B' \to M^{\ast} \to M[1].$$

Rewrite the last triangle as
$$M^{\ast} \overset{h}{\rightarrow} \tau M \overset{k}{\rightarrow} \tau B' \to M^{\ast}[1]$$
where we use that $M^{\ast}[1] = \tau^{-1}M^{\ast}$ and
$\tau M^{\ast} = M^{\ast}[1]$ in $\C$.
This gives rise to an exact sequence of functors
$$\Hom_{\C}(\ ,M^{\ast}) \to \Hom_{\C}(\ ,\tau M) \to \Hom_{\C}(\ ,\tau B') \to \Hom_{\C}(\ ,M^{\ast}[1]) 
\to \cdots.$$
Assume that $X$ is an indecomposable object which is not isomorphic to $M^{\ast}$,
and which is in $\Supp_{\C}(\ ,\tau M)$, and let $s \in \Hom_{\C}(X, \tau M)$ be a non-zero map.

If $ks \colon X \to \tau B'$ is not zero, then 
$X$ is in $\Supp_{\C}(\ , \tau B')$. If
$ks = 0$, then there is some $s' \colon X \to M^{\ast}$ such that 
$s = hs'$.
Denote by 
$$\tau M^{\ast} \overset{a}{\rightarrow} A \overset{r}{\rightarrow}M^{\ast} 
\to \tau M^{\ast}[1]$$
the almost split triangle in $\C$ for $M^{\ast}$.
Since $X \not \simeq M^{\ast}$, there is some map $s'' \colon X \to A$ such
that $s' = r{s''}$. Consider the commutative diagram in
Figure~\ref{commdiag2},
\begin{figure}[htbp]
\beginpicture
\setcoordinatesystem units <0.8cm,0.6cm>             
\setplotarea x from -6.2 to 5, y from 0 to 3       

\put{$\tau M^{\ast}$}[c] at 0 0
\put{$\tau B$}[c] at 2 0
\put{$\tau M$}[c] at 4 0
\put{$\tau M^{\ast}$}[c] at 0 2
\put{$A$}[c] at 2 2
\put{$M^{\ast}$} at 4 2
\scriptsize{
\put{$\tau g$}[c] at 0.8 0.3
\put{$a$}[c] at 0.8 2.3
\put{$\tau f$}[c] at 2.8 0.3
\put{$r$}[c] at 2.8 2.3
\put{$b_1$}[c] at 2.4 1
\put{$b_2$}[c] at 4.4 1
}
\arrow <1.8mm> [.5,1] from 0.35 0.03 to 1.3 0.03
\arrow <1.8mm> [.5,1] from 2.25 0.03 to 3.3 0.03
\arrow <1.8mm> [.5,1] from 0.35 2.03 to 1.3 2.03
\arrow <1.8mm> [.5,1] from 2.25 2.03 to 3.3 2.03

\setlinear \plot  -0.2 0.4  -0.2 1.6  / %
\setlinear \plot  -0.3 0.4  -0.3 1.6  / %

\setdashes <0.5mm,0.5mm>
\arrow <1.8mm> [.5,1] from 1.8 1.6 to 1.8 0.4
\arrow <1.8mm> [.5,1] from 3.7 1.6 to 3.7 0.4
\endpicture
\caption{Commutative diagram for the proof of Lemma~\ref{supp}}
\label{commdiag2}
\end{figure}
where the map $b_1 \colon A \to \tau B$ exists since the first triangle is almost split
and the second one is not split, and $b_2$ is then the induced map.

We claim that the map $b_1 {s''} \colon X \to A \to \tau B$ is nonzero.
Note that $(\tau f) b_1 {s''} = b_2 r {s''} = b_2 s'$.
Since $\Ext^1_{\C}(M,M^{\ast})$ has dimension one over $D_M$,
it follows that $\Hom_{\C}(M^{\ast}, \tau M)$ also has dimension one.
Since $b_2$ and $h$ are both nonzero elements in $\Hom_{\C}(M^{\ast}, \tau M)$,
it follows that there is a nonzero map $\phi \colon \tau M \to \tau M$, 
necessarily an isomorphism, such 
that $b_2 = \phi k$. Hence $b_2 s' = \phi h s' = \phi s \neq 0$,
and consequently $b_1 {s''} \neq 0$.
This finishes the proof of the lemma.
\end{proof}

We need some additional preliminary results.

\begin{lem} \label{morthogonal}
Let the assumptions and notation be as before.
\begin{itemize}
\item[(a)]$\Ext^1_{\C}(M,T_i) = 0$ for any indecomposable summand $T_i$ of $T$ 
which is not isomorphic to $M^{\ast}$.
\item[(b)]$\Ext^1_{\C}(M^{\ast},T_i) = 0$ for any indecomposable summand $T_i$ of $T$ 
which is not isomorphic to $M$.
\end{itemize}
\end{lem}

\begin{proof}
(a) Assume to the contrary that $\Ext^1_{\C}(M,T_1) \neq 0$ for some
$T_1$ an indecomposable summand of $T$, with $T_1 \not \simeq M^{\ast}$.
We have $\Ext^1_{\C}(M,T_1) \simeq 
D\Hom_{\C}(T_1, \tau M) \neq 0$, and hence by Lemma~\ref{supp}, either
$\Hom_{\C}(T_1, \tau B) \neq 0$ or $\Hom_{\C}(T_1, \tau B') \neq 0$, so that
$\Ext^1_{\C}(B,T_1) \neq 0$ or $\Ext^1_{\C}(B',T_1) \neq 0$.
But this contradicts the fact that $B \coprod B' \coprod T'$ is exceptional,
and the claim follows.

(b) The proof is dual to the proof of (a).
\end{proof}

We can now get the following.

\begin{lem}
If $M^{\ast}$ is not a direct summand of $T$, then $M$ is a direct
summand of $T$, and if $T = M^k \coprod \overline{T}$ (with $M$ not
a direct summand of $\overline{T}$), then $M^{\ast}\coprod \overline{T}$
is also a tilting object.
\end{lem}

\begin{proof}
Assume that $M$ and $M^{\ast}$ are not summands of $T$. Then
by Lemma~\ref{morthogonal}, $T \coprod M$ is exceptional, contradicting the
fact that $T$ is a tilting object.

Assume still that $M^{\ast}$ is not a summand of $T$, so that 
$T = M^k \coprod \overline{T}$ where $M$ is not a summand of $\overline{T}$
and $k>0$. By Lemma~\ref{morthogonal}, $M^{\ast}\coprod \overline{T}$ is an
exceptional object with the ``correct" number of indecomposable non-isomorphic
direct summands, and is hence a tilting object.  
\end{proof}

Summarising, we now have the following.

\begin{thm} \label{exchangepair}
Two exceptional indecomposable objects $M$ and $M^{\ast}$ form an exchange
pair if and only if
$\dim_{D_M}\Ext^1_{\C}(M,M^{\ast}) = 1 = \dim_{D_{M^{\ast}}}\Ext^1_{\C}(M^{\ast},M)$.
\end{thm}

The following example shows that it is necessary to assume that both
$\Ext^1$-spaces 
are one-dimensional, that is, one is not the consequence of the other.

\subsection*{Example}
\begin{figure}
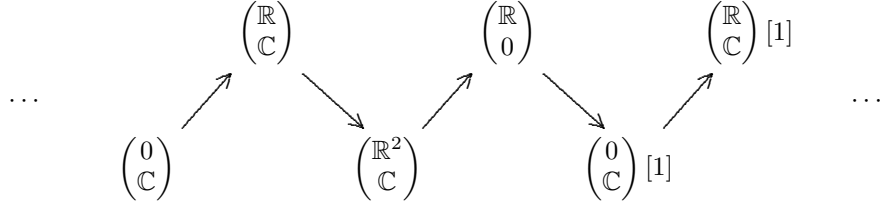

\beginpicture
\setcoordinatesystem units <0.8cm,0.9cm>             
\setplotarea x from -3.5 to 4, y from -1 to 5       

\put{$\begin{pmatrix} 0 \\ \mathbb{C} \end{pmatrix}$}[c] at 0 0
\put{$\begin{pmatrix} \mathbb{R} \\ \mathbb{C} \end{pmatrix}$}[c] at 2 2
\put{$\begin{pmatrix} \mathbb{R}^2 \\ \mathbb{C} \end{pmatrix}$}[c] at 4 0
\put{$\begin{pmatrix} \mathbb{R} \\ 0 \end{pmatrix}$}[c] at 6 2
\put{$\begin{pmatrix} 0 \\ \mathbb{C} \end{pmatrix}[1]$}[c] at 8 0
\put{$\begin{pmatrix} \mathbb{R} \\ \mathbb{C} \end{pmatrix}[1]$}[c] at 10 2

\put{$\cdots$} at -2 1
\put{$\cdots$} at 12 1

\arrow <1.8mm> [.5,1] from 0.6 0.6 to 1.4 1.4
\arrow <1.8mm> [.5,1] from 4.6 0.6 to 5.4 1.4
\arrow <1.8mm> [.5,1] from 8.6 0.6 to 9.4 1.4
\arrow <1.8mm> [.5,1] from 2.6 1.4 to 3.6 0.6
\arrow <1.8mm> [.5,1] from 6.6 1.4 to 7.6 0.6

\endpicture
\caption{The AR-quiver of $D^b(H)$}
\label{nonsymmetricexample}
\end{figure}

Consider the ring 
$$ H =\left( \begin{array}{cc}
\mathbb{R} & 0 \\
_{\mathbb{R}}{\mathbb{C}}_{\mathbb{R}} & \mathbb{C}  
\end{array} \right). 
$$
The AR-quiver of $D^b(H)$ is shown in
Figure~\ref{nonsymmetricexample}.
We have
$$\dim_{\mathbb{C}} \Ext^1_{\C}(\begin{pmatrix} \mathbb{R} \\ 0 \end{pmatrix} ,
\begin{pmatrix} 0 \\ \mathbb{C} \end{pmatrix}) =1$$
and
$$\dim_{\mathbb{R}} \Ext^1_{\C}( \begin{pmatrix} \mathbb{R} \\ 0 \end{pmatrix},
 \begin{pmatrix} 0 \\ \mathbb{C} \end{pmatrix}) =2.$$

Hence $\{ \begin{pmatrix} \mathbb{R} \\ 0 \end{pmatrix},
 \begin{pmatrix} 0 \\ \mathbb{C} \end{pmatrix} \}$
is not an exchange pair.

Finally, suppose that $H$ is the path algebra of a quiver of simply-laced
Dynkin type $\Delta$ with an alternating orientation. Let
$\mathcal{A}(\Delta)$ denote the corresponding
cluster algebra. By Proposition~\ref{objectbijection}, we know that there is
a 1--1 correspondence between the cluster variables of
$\mathcal{A}(\Delta)$ and $\ind\C$. By Theorem~\ref{clusterbijection} we know
that this induces a bijection between the basic tilting objects of $\C$ and
the clusters of $\mathcal{A}(\Delta)$. We have the following interpretation of
Theorem~\ref{exchangepair}.

\begin{theorem} \cite[3.5,4.4]{fz2}
Suppose $\mathcal{A}(\Delta)$ is the cluster algebra associated to an
arbitrary Dynkin diagram of simply-laced type.
Let $x,y$ be two cluster variables of $\mathcal{A}(\Delta)$. Then
$x,y$ form an exchange pair
if and only if their compatibility degree is equal to $1$.
\end{theorem}

\section{Graphical calculus} \label{graphicalcalculus}

In this section, we assume the quiver $Q$ to be of simply-laced Dynkin type.
We shall give a graphical calculus for computing the triangles in
Section~\ref{approx} (see Theorem~\ref{magictriangle} and the comment
afterwards).

Suppose that $M,M^{\ast}$ are indecomposable objects of $\C=D^b(kQ)/F$.
We know that $\End_{\C}(M)\simeq \End_{\C}(M')\simeq k$ --- use Proposition~\ref{symm}(c) and
the fact that every indecomposable object in $\C$ is in the $\tau$-orbit of
an indecomposable projective module.
Suppose that $\Ext^1_{\C}(M^{\ast},M)$ is one-dimensional over $k$.
We know by Theorem~\ref{exchangepair} that this is the equivalent to assuming that
$M,M^{\ast}$ are the two complements of an almost complete basic tilting object $\overline{T}$
of $\C$. We would like to construct triangles:
\begin{equation} \label{magic1}
M^{\ast} \overset{g}{\rightarrow} B \overset{f}{\rightarrow} M \to M^{\ast}[1]
\end{equation}
and
\begin{equation}\label{magic2}
M \overset{u}{\rightarrow} B' \overset{v}{\rightarrow} M^{\ast} \to M[1]
\end{equation}
where $f:B\rightarrow M$ is a minimal right $\add\overline{T}$-approximation of $M$,
and $u:M\rightarrow B'$ is a minimal left $\add\overline{T}$-approximation of $M$.

Without loss of generality (by applying APR-tilts if necessary), we can assume
that $M$ is a simple projective $kQ$-module $P$. Suppose first that
$M^{\ast}=P'[1]$ is the shift of an indecomposable projective
$kQ$-module $P'$.
Then 
$$\Ext^1_{\C}(M^{\ast},M)=\Hom_{kQ}(P',P)$$
(see Proposition~\ref{symm}(d)),
since $P'$ is projective. But, since $P$ is simple projective, this is
non-zero (and necessarily one-dimensional) if and only if $P\simeq P'$.
By applying the autoequivalence $\tau_{\C}^{-1}$ to $P$ and $P'[1]$
we are reduced to the situation where $M$ and $M^{\ast}$ are the start
and end terms respectively of an almost split sequence of $kQ$-modules. We are
then in the case 
discussed after the proof of Proposition~\ref{dimensionone}, and we see that
$B=0$ and $B'$ is the middle term of the almost split sequence involving
$M$ and $M^{\ast}$.

We are now left with the case where $M$ and $M^{\ast}$ are both modules
over $kQ$, with $M$ projective. By Proposition~\ref{symm}(c), we have:
$$\Ext^1_{\C}(M^{\ast},M)\simeq \Ext^1_{kQ}(M^{\ast},M)\coprod \Ext^1_{kQ}(M,M^{\ast}).$$
Since $M$ is projective, $\Ext^1_{kQ}(M,M^{\ast})=0$.
Then we have a unique non-trivial extension
$$0\rightarrow M\rightarrow E\rightarrow M^{\ast}\rightarrow 0$$
of $kQ$-modules. There is a corresponding triangle
$$M\rightarrow E\rightarrow M^{\ast}\rightarrow M[1]$$
in $\D$ which induces a non-split triangle in $\C$. Since the
triangle~\eqref{magic2} is (up to isomorphism) the unique non-split triangle in
$\C$ with start term $M$ and end term $M^{\ast}$, we have that $E$ is
isomorphic to $B'$. In the case where
$\Ext^1_{kQ}(M,M^{\ast})\simeq k$ we obtain a middle term isomorphic to $B$.
Note that by switching the roles of $M$ and $M^{\ast}$ (using
Proposition~\ref{symm}(b)), and applying $\tau_{\C}$ as appropriate, we can
compute both middle terms $B$ and $B'$ using
the module category alone. We have reduced the problem to the following:

\begin{problem} \label{sesproblem}
Let $Q$ be a simply-laced Dynkin quiver, and let $M,M^{\ast}$ be indecomposable
$kQ$-modules satisfying $\Ext^1_{kQ}(M,M^{\ast})\simeq k$ (and therefore
$\Ext^1_{kQ}(M^{\ast},M)=0$). Compute the middle term of the unique non-trivial
extension represented by a non-zero element of $\Ext^1_{kQ}(M,M^{\ast})$.
\end{problem}

Let $M$ and $M^{\ast}$ be indecomposable $kQ$-modules such
that $\dim\Ext_{kQ}^1(M,M^{\ast})=1$, and let
$$\zeta\, :\, 0\rightarrow M^{\ast}\rightarrow X\rightarrow M\rightarrow 0$$
be the unique non-trivial extension mentioned above.
We will now develop a graphical method (in terms of the AR-quiver) for the
determination of $X$.
We recall that the {\em starting function} $s_U$ of an indecomposable $kQ$-module
$U$ is defined as the function $V\mapsto\dim \Hom_{kQ}(U,V)$ on
indecomposable $kQ$-modules. All such starting functions are depicted
in~\cite{bo2}. Similarly, the ending function $e_U$ is defined as the function
$V\mapsto\dim \Hom_{kQ}(V,U)$.

\begin{lemma} \label{smalldim}
Let $U$ and $V$ be indecomposable representations of $kQ$ such
that $\Hom_{kQ}(U,V)\not=0$ and $\Ext_{kQ}^1(V,U)=0$. Then
$\dim\Hom_{kQ}(U,V)=1$.
\end{lemma}

\begin{proof}
The above condition translates to $$s_U(V)\not=0,\;\; s_U(\tau V)=0$$
by the AR-formula. Now direct inspection of the tables
in~\cite{bo2} gives the above result.
\end{proof}

This result can also be established in a theoretical way, using the
result~\cite[Corollary to main Theorem]{vh}.
Since the table in~\cite{bo2} will play a central role in the following, the
above proof is more adapted to the theme of this section.

\begin{prop} Let $M$, $M^{\ast}$ and $X$ be as above. Then $X$ is the direct sum of one copy of each indecomposable $kQ$-module $V$ fulfilling
$$\Hom_{kQ}(M^{\ast},V)\not=0\not=\Hom_{kQ}(V,M)\mbox{ and }\Ext_{kQ}^1(V,M^{\ast})=0=\Ext_{kQ}^1(M,V).$$
\end{prop}

\begin{proof}
Let $V$ be an indecomposable direct summand of $X$. We first show that the 
stated homological conditions on $V$ are satisfied. If $\Hom_{kQ}(M^{\ast},V)=0$, 
then $V$ has to appear as a direct summand of $M$. Since $M$ is indecomposable, 
this implies $M=V$, and the sequence $\zeta$ splits, a contradiction. 
Thus $\Hom_{kQ}(M^{\ast},V)\not=0$.\\[1ex]
In the induced exact sequence
$$\Hom_{kQ}(M^{\ast},M^{\ast})\stackrel{d}{\rightarrow}\Ext_{kQ}^1(M,M^{\ast})\rightarrow
\Ext_{kQ}^1(X,M^{\ast})\rightarrow\Ext_{kQ}^1(M^{\ast},M^{\ast})=0,$$
the map $d$ is surjective, since $\Ext_{kQ}^1(M,M^{\ast})$ is one-dimensional and 
the sequence $\zeta$ is non-split. Thus $\Ext_{kQ}^1(X,M^{\ast})=0$, and in particular 
$\Ext_{kQ}^1(V,M^{\ast})=0$.  
We can argue dually to obtain the other two conditions on $V$.\\[1ex]
Enumerate the isomorphism classes of indecomposables with the above 
properties as $\{V_1,\ldots,V_s\}$; thus we can write $X=\oplus_{i=1}^sV_i^{m_i}$, 
and we have to prove that $m_i=1$ for all $i=1\ldots s$. Consider the induced exact 
sequence
$$0\rightarrow\Hom_{kQ}(M,X)\rightarrow\Hom_{kQ}(X,X)\rightarrow
\Hom_{kQ}(M^{\ast},X)\rightarrow\Ext_{kQ}^1(M,X).$$
From the above, we can conclude that $\Ext^1_{kQ}(M,X)=0$. 
Since any $V_i$ maps to $M$, we also have $\Hom_{kQ}(M,X)=0$, thus $\Hom_{kQ}(M,V_i)=0$ 
since the category $\bmod kQ$ is representation-directed. We arrive at an 
isomorphism $$\Hom_{kQ}(X,X)\simeq\Hom_{kQ}(M^{\ast},X).$$
Since ${\rm End}_{kQ}(X)$ contains the semisimple ring 
$\oplus_{i=1}^sM_{m_i}({\rm End}_{kQ}(V_i))$ as a subring, we can estimate:
$$\sum_{i=1}^sm_i^2\leq\sum_{i,j=1}^sm_im_j\dim\Hom_{kQ}(V_i,V_j)=\dim\Hom_{kQ}(X,X)=$$
$$=\dim\Hom_{kQ}(M^{\ast},X)=\sum_{i=1}^sm_i\dim\Hom_{kQ}(M^{\ast},V_i)=\sum_{i=1}^sm_i,$$
using Lemma~\ref{smalldim}.
Thus $m_i\in\{0,1\}$ for all $i=1\ldots s$, and $\Hom_{kQ}(V_i,V_j)=0$ whenever 
$i\not=j$ and $m_i=1=m_j$.\\[1ex]
Similarly, we see that for each $i=1\ldots s$, we have isomorphisms
$$\Hom_{kQ}(X,V_i)\simeq\Hom_{kQ}(M^{\ast},V_i)\mbox{ and }\Hom_{kQ}(V_i,X)\simeq
\Hom_{kQ}(V_i,M).$$
Given a fixed $V_i$, we choose non-zero maps 
$f:M^{\ast}\rightarrow V_i$ and $g:V_i\rightarrow M$. The above isomorphisms yield 
factorisations $f=r\alpha$ and $g=\beta s$, where $\alpha:M^{\ast}\rightarrow X$
and $\beta:X\rightarrow M$ are the maps in the short exact sequence $\zeta$.
Since $r\not=0$ and $s\not=0$ we can choose summands $V_j$ and $V_k$ of $X$
such that $r_j\alpha\not=0$ and $\beta s_k\not=0$, where $r_j$ is the
restriction of $r$ to $V_j$ and $s_k$ is the composition of $s$ with the
projection onto $V_k$.

It is enough to prove that $s_kr_j\not=0$. Then, since $V_j$ and $V_k$ are
direct summands of $X$, we obtain from the above that $j=k$, and therefore
that $i=j=k$ since there are no oriented cycles of homomorphisms in
the category of $kQ$-modules.

We have $\Ext_{kQ}^1(M,V_i)=0$ by assumption, and $\Ext_{kQ}^1(M^{\ast},V_i)=0$ since $kQ$ is
representation-directed and $\Hom_{kQ}(M^{\ast},V_i)\not=0$. These two facts together imply 
$\Ext_{kQ}^1(X,V_i)=0$, 
thus in particular $\Ext_{kQ}^1(V_k,V_i)=0$, since $V_k$ is a direct 
summand of $X$. This vanishing condition allows us to apply the Happel-Ringel 
Lemma~\cite{hr} to conclude that $s_k\not=0$ must be mono or epi. If $s_k$ is
mono, then $s_kr_j\not=0$ since $r_j\not=0$, and we are done. So assume that
$s_k$ is epi. By possibly applying the AR-translate, we can assume without
loss of generality that $M^{\ast}$ is projective. This provides us with a
surjection
$$(s_k\circ\_):\Hom_{kQ}(M^{\ast},V_i)\rightarrow\Hom_{kQ}(M^{\ast},V_k),$$
thus an isomorphism since both spaces are one-dimensional by
Lemma~\ref{smalldim}. But this implies that $s_kr_j\not=0$.
This finishes the proof.
\end{proof}

The starting and ending functions of an indecomposable $kQ$-module $U$ can be
computed in terms of the AR-quiver: the function $s_U$ is determined by defining
$s_U(V)=1$ on the slice starting in $U$, and by additivity
$s_U(\tau^{-1}(V))=\sum_is_U(C_i)-s_U(V)$ for a mesh
$V\rightarrow\oplus_iC_i\rightarrow\tau^{-1}V$.
We can now define:

\begin{definition} The starting frame $F_s(U)$ (resp.~the ending frame $F_e(U)$)
of an indecomposable $kQ$-module $U$ consists of all vertices $V$ of the AR-quiver
such that $s_U(V)\not=0=s_U(\tau V)$ (resp.~$e_U(V)\not=0=e_U(\tau^{-1}V)$).
\end{definition}

As an immediate corollary to the above proposition, we get:

\begin{cor} Given indecomposables $M$ and $M^{\ast}$ such that $\Ext_{kQ}^1(M,M^{\ast})$ is 
one-dimensional, the unique non-trivial extension $X$ of $M$ by $M^{\ast}$ is given as 
the direct sum of all indecomposables belonging to the intersection $F_s(M^{\ast})\cap F_e(M)$.
\end{cor}

The starting and ending frames can now be worked out using the tables in \cite{bo2}.
In type $A$, they are easily seen to coincide with the slice starting (resp.~ending) in 
an indecomposable. For type $D$ and $E$, the frames look in general more complicated. 
Below, we first show two "typical" examples in type $D_8$. The starting frame of the 
respective minimal vertex of the picture is shown (the solid circles), embedded in a 
portion of the AR-quiver.

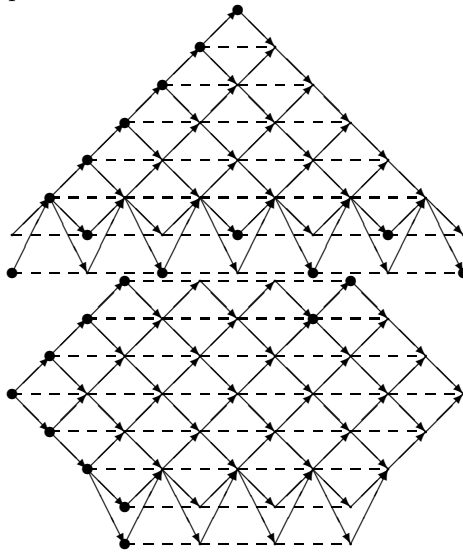
\setlength{\unitlength}{0.5cm}
\begin{center}
\begin{picture}(12,7)

\multiput(0,0)(2,0){7}{\circle*{0.1}}
\multiput(0,1)(2,0){7}{\circle*{0.1}}
\multiput(1,2)(2,0){6}{\circle*{0.1}}
\multiput(2,3)(2,0){5}{\circle*{0.1}}
\multiput(3,4)(2,0){4}{\circle*{0.1}}
\multiput(4,5)(2,0){3}{\circle*{0.1}}
\multiput(5,6)(2,0){2}{\circle*{0.1}}
\put(6,7){\circle*{0.1}}

\multiput(0,0)(2,0){6}{\vector(1,2){1}}
\multiput(0,1)(2,0){6}{\vector(1,1){1}}
\multiput(1,2)(2,0){5}{\vector(1,1){1}}
\multiput(1,2)(2,0){6}{\vector(1,-1){1}}
\multiput(1,2)(2,0){6}{\vector(1,-2){1}}
\multiput(2,3)(2,0){4}{\vector(1,1){1}}
\multiput(2,3)(2,0){5}{\vector(1,-1){1}}
\multiput(3,4)(2,0){3}{\vector(1,1){1}}
\multiput(3,4)(2,0){4}{\vector(1,-1){1}}
\multiput(4,5)(2,0){2}{\vector(1,1){1}}
\multiput(4,5)(2,0){3}{\vector(1,-1){1}}
\put(5,6){\vector(1,1){1}}
\multiput(5,6)(2,0){2}{\vector(1,-1){1}}
\put(6,7){\vector(1,-1){1}}

\multiput(0,0)(0.5,0){24}{\line(1,0){0.25}}
\multiput(0,1)(0.5,0){24}{\line(1,0){0.25}}
\multiput(1,2)(0.5,0){20}{\line(1,0){0.25}}
\multiput(2,3)(0.5,0){16}{\line(1,0){0.25}}
\multiput(3,4)(0.5,0){12}{\line(1,0){0.25}}
\multiput(4,5)(0.5,0){8}{\line(1,0){0.25}}
\multiput(5,6)(0.5,0){4}{\line(1,0){0.25}}

\put(0,0){\circle*{0.3}}
\put(1,2){\circle*{0.3}}
\put(2,1){\circle*{0.3}}
\put(2,3){\circle*{0.3}}
\put(3,4){\circle*{0.3}}
\put(4,0){\circle*{0.3}}
\put(4,5){\circle*{0.3}}
\put(5,6){\circle*{0.3}}
\put(6,1){\circle*{0.3}}
\put(6,7){\circle*{0.3}}
\put(8,0){\circle*{0.3}}
\put(10,1){\circle*{0.3}}
\put(12,0){\circle*{0.3}}

\end{picture}
\end{center}

\setlength{\unitlength}{0.5cm}
\begin{center}
\begin{picture}(12,7)

\multiput(3,0)(2,0){4}{\circle*{0.1}}
\multiput(3,1)(2,0){4}{\circle*{0.1}}
\multiput(2,2)(2,0){5}{\circle*{0.1}}
\multiput(1,3)(2,0){6}{\circle*{0.1}}
\multiput(0,4)(2,0){7}{\circle*{0.1}}
\multiput(1,5)(2,0){6}{\circle*{0.1}}
\multiput(2,6)(2,0){5}{\circle*{0.1}}
\multiput(3,7)(2,0){4}{\circle*{0.1}}

\multiput(3,0)(2,0){4}{\vector(1,2){1}}
\multiput(3,1)(2,0){4}{\vector(1,1){1}}
\multiput(2,2)(2,0){5}{\vector(1,1){1}}
\multiput(2,2)(2,0){4}{\vector(1,-1){1}}
\multiput(2,2)(2,0){4}{\vector(1,-2){1}}
\multiput(1,3)(2,0){6}{\vector(1,1){1}}
\multiput(1,3)(2,0){5}{\vector(1,-1){1}}
\multiput(0,4)(2,0){6}{\vector(1,1){1}}
\multiput(0,4)(2,0){6}{\vector(1,-1){1}}
\multiput(1,5)(2,0){5}{\vector(1,1){1}}
\multiput(1,5)(2,0){6}{\vector(1,-1){1}}
\multiput(2,6)(2,0){4}{\vector(1,1){1}}
\multiput(2,6)(2,0){5}{\vector(1,-1){1}}
\multiput(3,7)(2,0){4}{\vector(1,-1){1}}

\multiput(3,0)(0.5,0){12}{\line(1,0){0.25}}
\multiput(3,1)(0.5,0){12}{\line(1,0){0.25}}
\multiput(2,2)(0.5,0){16}{\line(1,0){0.25}}
\multiput(1,3)(0.5,0){20}{\line(1,0){0.25}}
\multiput(0,4)(0.5,0){24}{\line(1,0){0.25}}
\multiput(1,5)(0.5,0){20}{\line(1,0){0.25}}
\multiput(2,6)(0.5,0){16}{\line(1,0){0.25}}
\multiput(3,7)(0.5,0){12}{\line(1,0){0.25}}

\put(3,0){\circle*{0.3}}
\put(3,1){\circle*{0.3}}
\put(2,2){\circle*{0.3}}
\put(1,3){\circle*{0.3}}
\put(0,4){\circle*{0.3}}
\put(1,5){\circle*{0.3}}
\put(2,6){\circle*{0.3}}
\put(3,7){\circle*{0.3}}
\put(8,6){\circle*{0.3}}
\put(9,7){\circle*{0.3}}

\end{picture}
\end{center}

Finally, we show the "most complicated" starting frame in type $E_8$:

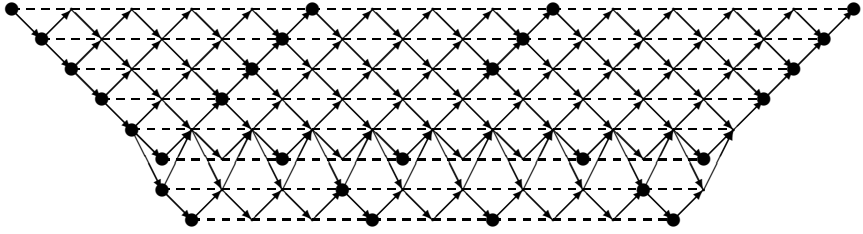
\setlength{\unitlength}{0.4cm}
\begin{center}
\begin{picture}(28,7)

\multiput(6,0)(2,0){9}{\circle*{0.1}}
\multiput(5,1)(2,0){10}{\circle*{0.1}}
\multiput(5,2)(2,0){10}{\circle*{0.1}}
\multiput(4,3)(2,0){11}{\circle*{0.1}}
\multiput(3,4)(2,0){12}{\circle*{0.1}}
\multiput(2,5)(2,0){13}{\circle*{0.1}}
\multiput(1,6)(2,0){14}{\circle*{0.1}}
\multiput(0,7)(2,0){15}{\circle*{0.1}}

\multiput(6,0)(2,0){9}{\vector(1,1){1}}
\multiput(5,1)(2,0){10}{\vector(1,2){1}}
\multiput(5,1)(2,0){9}{\vector(1,-1){1}}
\multiput(5,2)(2,0){10}{\vector(1,1){1}}
\multiput(4,3)(2,0){10}{\vector(1,-2){1}}
\multiput(4,3)(2,0){10}{\vector(1,-1){1}}
\multiput(4,3)(2,0){11}{\vector(1,1){1}}
\multiput(3,4)(2,0){11}{\vector(1,-1){1}}
\multiput(3,4)(2,0){12}{\vector(1,1){1}}
\multiput(2,5)(2,0){12}{\vector(1,-1){1}}
\multiput(2,5)(2,0){13}{\vector(1,1){1}}
\multiput(1,6)(2,0){13}{\vector(1,-1){1}}
\multiput(1,6)(2,0){14}{\vector(1,1){1}}
\multiput(0,7)(2,0){14}{\vector(1,-1){1}}

\multiput(6,0)(0.5,0){32}{\line(1,0){0.25}}
\multiput(5,1)(0.5,0){36}{\line(1,0){0.25}}
\multiput(5,2)(0.5,0){36}{\line(1,0){0.25}}
\multiput(4,3)(0.5,0){40}{\line(1,0){0.25}}
\multiput(3,4)(0.5,0){44}{\line(1,0){0.25}}
\multiput(2,5)(0.5,0){48}{\line(1,0){0.25}}
\multiput(1,6)(0.5,0){52}{\line(1,0){0.25}}
\multiput(0,7)(0.5,0){56}{\line(1,0){0.25}}

\put(0,7){\circle*{0.4}}
\put(1,6){\circle*{0.4}}
\put(2,5){\circle*{0.4}}
\put(3,4){\circle*{0.4}}
\put(4,3){\circle*{0.4}}
\put(5,2){\circle*{0.4}}
\put(5,1){\circle*{0.4}}
\put(6,0){\circle*{0.4}}
\put(7,4){\circle*{0.4}}
\put(8,5){\circle*{0.4}}
\put(9,6){\circle*{0.4}}
\put(9,2){\circle*{0.4}}
\put(10,7){\circle*{0.4}}
\put(11,1){\circle*{0.4}}
\put(12,0){\circle*{0.4}}
\put(13,2){\circle*{0.4}}
\put(16,5){\circle*{0.4}}
\put(16,0){\circle*{0.4}}
\put(17,6){\circle*{0.4}}
\put(18,7){\circle*{0.4}}
\put(19,2){\circle*{0.4}}
\put(21,1){\circle*{0.4}}
\put(22,0){\circle*{0.4}}
\put(23,2){\circle*{0.4}}
\put(25,4){\circle*{0.4}}
\put(26,5){\circle*{0.4}}
\put(27,6){\circle*{0.4}}
\put(28,7){\circle*{0.4}}

\end{picture}
\end{center}

\section{Interpretation and conjectures}

In this section we will consider further links with cluster algebras,
including interpretations of some of the preceding results. We will
make some conjectures in this direction and also provide some examples
giving supporting evidence for the conjectures.

Let $H$ be a finite dimensional hereditary algebra, with quiver $\Gamma$.
For vertices $i$ and $j$ of $\Gamma$, let $n_{ij}$ denote the number of
arrows from $i$ to $j$ in $\Gamma$.
Let $X$ be the matrix with rows and columns indexed by the
vertices of $\Gamma$ (we choose a total ordering):
$$x_{ij}=\left\{
\begin{array}{ll}
n_{ij} & n_{ij}\not=0, \\
-n_{ji} & n_{ij}=0.
\end{array}\right.
$$
Let $\A(H)$ be the corresponding cluster algebra.
Let $\C$ be the cluster category associated to $H$.

\begin{conj} \label{correspondenceconjecture}
There is a 1--1 correspondence between the cluster variables of
$\A(H)$ and $\ind \C$ inducing a 1--1 correspondence
between the clusters of $\A(H)$ and the basic tilting objects in $\C$.
\end{conj}

We have seen (see Section~\ref{connections}) that this conjecture holds in the
case where $H$ is the path algebra of a simply-laced Dynkin quiver. In this
case, we make a further conjecture:

\begin{conj} \label{linkwithcss}
Let $C$ be a cluster of the cluster algebra of simply-laced Dynkin type,
and let $T$ be the corresponding tilting object of the cluster category
$\C$ of the same type. Let $A_C$ denote the algebra associated to $C$
in~\cite[Section 1]{ccs} (its module category is denoted $\Mod Q_C$ there).
Then $\End_{\C}(T)^{\op}$ is isomorphic to $A_C$.
\end{conj}

Suppose that $\overline{T}$ is an almost complete basic tilting object of $\C$.
Let $M,M^{\ast}$ be the complements of $\overline{T}$, and let
$$M^{\ast} \overset{g}{\rightarrow} B \overset{f}{\rightarrow} M \to M^{\ast}[1]$$
and
$$M \overset{u}{\rightarrow} B' \overset{v}{\rightarrow} M^{\ast\ast} \to M[1]$$
be the triangles~\eqref{magictriangle1} and~\eqref{magictriangle2} from
Section~\ref{approx}.

We make the following conjecture:
\begin{conj} \label{exchangeconjecture}
In the above situation, let $B=\coprod_{i\in I}B_i^{d_i}$ (respectively,
$B'=\coprod_{j\in J}(B'_j)^{e_j}$) be the direct sum decomposition of $B$
(respectively, $B'$), where the $B_i$ are all non-isomorphic and the $B'_j$ are
all non-isomorphic. Let $x,x'$ be the cluster variables corresponding to
$M,M^{\ast}$, and for $i\in I$ (respectively, $j\in J$) let $x_i$ (respectively,
$x'_j$) be the cluster variable corresponding to $B_i$ (respectively,
$B'_j$). Then the exchange relation in the cluster algebra $\A(H)$
(see equation~\eqref{exchangerelation} in the introduction)
takes the form:
$$xx'=\prod_{i\in I}x_i^{d_i}+\prod_{i\in I'}(x'_j)^{e_j}.$$
In particular, $B$ and $B'$ should have no common direct summands.
\end{conj}

We note that in the simply-laced Dynkin case, this conjecture can be reformulated,
via the discussion in Section~\ref{graphicalcalculus}, to give a conjecture
providing a direct interpretation of the cluster exchange relation in terms of
short exact sequences of $kQ$-modules (see Problem~\ref{sesproblem}).
We also note that if Conjecture~\ref{exchangeconjecture} holds then it can be
seen that the rule for matrix
mutation (see the introduction) describes the change in the quiver of
the algebra $\End_{\C}(T)^{\op}$ when one indecomposable direct summand of the
basic tilting object $T$ is exchanged for another. We give an example of this
below. Finally, we remark that if Conjectures~\ref{correspondenceconjecture}
and~\ref{exchangeconjecture} both hold, then
Theorems~1.11 and~1.12 (without coefficients) in~\cite{fz2} hold for the
corresponding cluster algebra.

\subsection*{Example}
\begin{figure}
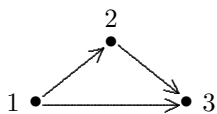

\beginpicture
\setcoordinatesystem units <0.5cm,0.4cm>             
\setplotarea x from -10.5 to 4, y from -1 to 5       

\put{$\bullet$}[c] at 0 0
\put{$\bullet$}[c] at 2 2
\put{$\bullet$}[c] at 4 0
\put{$1$}[c] at -0.6 0
\put{$2$}[c] at 2 2.8
\put{$3$}[c] at 4.6 0
\arrow <1.8mm> [.5,1] from 0.2 0.2 to 1.8 1.8
\arrow <1.8mm> [.5,1] from 2.2 1.9 to 3.8 0.3
\arrow <1.8mm> [.5,1] from 0.2 -0.1 to 3.8 -0.1

\endpicture
\caption{A quiver of type $\widetilde{A_2}$}
\label{a2tilde}
\end{figure}

Let $H$ be the path algebra of the quiver as shown in Figure~\ref{a2tilde}.
Then the corresponding cluster algebra $\A (H)$ has seed given by the
transcendence basis $\{u_1,u_2,u_3\}$ of
$\mathbb{Q}(u_1,u_2,u_3)$ and matrix
$$X=\left( \begin{array}{ccc}
0 & 1 & 1 \\
-1 & 0 & 1 \\
-1 & -1 & 0
\end{array}\right) .
$$
The corresponding cluster algebra was investigated in~\cite[7.8]{fz1} --- the 
{\em brick wall} example. Let $P_1,P_2,P_3$ denote the indecomposable
projective modules corresponding to the vertices of the graph of $H$.
Let $R$ denote the regular indecomposable module with dimension vector
$(1,0,1)$.
Then $T=\widetilde{P_1}\coprod \widetilde{P_2}\coprod \widetilde{P_3}$ is a
basic tilting object
of $\C$. Choosing $M=\widetilde{P_2}$ and $\overline{T}=\widetilde{P_1}\coprod \widetilde{P_3}$
we see that $\widetilde{P_2}\rightarrow \widetilde{P_1}$ is a minimal left 
$\add(\overline{T})$-approximation of $\widetilde{P_2}$ and obtain the triangle:
$$\widetilde{P_2}\rightarrow \widetilde{P_1}\rightarrow \widetilde{R}\rightarrow \widetilde{P_2}[1]$$
in $\C$. It follows that $T'=\widetilde{P_1}\coprod \widetilde{R}\coprod \widetilde{P_3}$ is
again a basic tilting object of $C$. The matrix $X'$ of the quiver of
$\End(T')^{\op}$ corresponding to $T'$ is:
$$X'=\left( \begin{array}{ccc}
0 & -1 & 2 \\
1 & 0 & -1 \\
-2 & 1 & 0
\end{array}\right) ,
$$
which is easily seen to be the mutation of the matrix $X$ at $2$.

Suppose that $H$ is a finite dimensional hereditary algebra, $\overline{T}$
is an almost complete basic tilting object of $\C$, and $M$ and $M^{\ast}$ are
the two complements of $\overline{T}$. Let $T = \overline{T} \coprod M$ and
$T' = \overline{T} \coprod M^{\ast}$ be the two completions of $\overline{T}$ to a
basic tilting object. Let $\G = \End_{\C}(\overline{T} \coprod M)^{\op}$ and 
$\G' = \End_{\C}(\overline{T} \coprod M^{\ast})^{\op}$ be the endomorphism algebras,
taken over $\C$.
Denote by $S_M$ (respectively, $S_{M^{\ast}}$) the simple top of the
$\G$-module $\Hom_{\C}(T,M)$ (respectively, the $\G'$-module
$\Hom_{\C}(T, M^{\ast})$).
Then we conjecture that the category of $\G$-modules and the category of
$\G'$-modules are related in the following way:

\begin{conj} \label{generalisedapr}
The categories $\mod \G / \add S_M$ and $\mod \G' / \add S_{M^{\ast}}$
are equivalent.
\end{conj}

This can be viewed as a generalisation of APR-tilting~\cite{apr}.

\subsection*{Example} \label{mutation}
\begin{figure}[htbp]
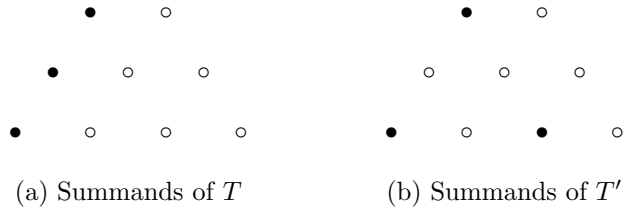

\beginpicture
\setcoordinatesystem units <0.5cm,0.4cm>             
\setplotarea x from -3 to 10, y from -10.5 to -1.5       

\put{$\bullet$}[c] at 1.8 -8
\put{$\circ$}[c] at 3.8 -8
\put{$\circ$}[c] at 5.8 -8
\put{$\circ$}[c] at 7.8 -8
\put{$\bullet$}[c] at 2.8 -6
\put{$\circ$}[c] at 4.8 -6
\put{$\circ$}[c] at 6.8 -6
\put{$\bullet$}[c] at 3.8 -4
\put{$\circ$}[c] at 5.8 -4

\put{(a) Summands of $T$}[c] at 4.8 -10

\setcoordinatesystem point at -10 0  

\put{$\bullet$}[c] at 1.8 -8
\put{$\circ$}[c] at 3.8 -8
\put{$\bullet$}[c] at 5.8 -8
\put{$\circ$}[c] at 7.8 -8
\put{$\circ$}[c] at 2.8 -6
\put{$\circ$}[c] at 4.8 -6
\put{$\circ$}[c] at 6.8 -6
\put{$\bullet$}[c] at 3.8 -4
\put{$\circ$}[c] at 5.8 -4

\put{(b) Summands of $T'$}[c] at 4.8 -10

\endpicture
\caption{Two basic tilting objects of $\C$ in type $A_3$}
\label{tiltexample}
\end{figure}

We give an example illustrating Conjecture~\ref{generalisedapr}. Take
$\Delta$ to be the Dynkin diagram of type $A_3$. Then the AR-quiver of
$\C$ is given in Figure~\ref{Cquiver}. Let $T$ be the direct sum
of the indecomposable objects corresponding to the filled-in circles
in Figure~\ref{tiltexample}(a) and let $T'$ be the direct sum of the
indecomposable objects corresponding to the filled-in circles in
Figure~\ref{tiltexample}(b). Thus $\overline{T}$ is the almost complete basic
tilting object which is the direct sum of the objects corresponding to the two
filled-in circles common to $T$ and $T'$.
Here we display the AR-quiver of $\C$ slightly differently in order
to demonstrate this example (noting that it appears on a M\"{o}bius band).

\begin{figure}[htbp]
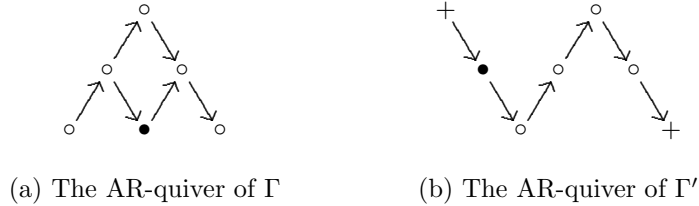


\beginpicture

\setcoordinatesystem units <0.5cm,0.4cm>             
\setplotarea x from -3 to 10, y from -10.5 to -1.5       

\put{$\circ$}[c] at 1.8 -8
\put{$\bullet$}[c] at 3.8 -8
\put{$\circ$}[c] at 5.8 -8
\put{$\circ$}[c] at 2.8 -6
\put{$\circ$}[c] at 4.8 -6
\put{$\circ$}[c] at 3.8 -4
\put{(a) The AR-quiver of $\Gamma$} at 3.8 -10

\setlinear \plot  2 -7.65   2.6 -6.4  / %
\setlinear \plot  2.6 -6.4 2.7 -6.7  / %
\setlinear \plot  2.6 -6.4 2.3 -6.5  / %

\setlinear \plot  4 -7.65   4.6 -6.4  / %
\setlinear \plot  4.6 -6.4 4.7 -6.7  / %
\setlinear \plot  4.6 -6.4 4.3 -6.5  / %

\setlinear \plot  3.6 -7.65 3.0 -6.4  / %
\setlinear \plot  3.6 -7.65 3.3 -7.55 / %
\setlinear \plot  3.6 -7.65 3.7 -7.35  / %

\setlinear \plot  5.6 -7.65 5.0 -6.4  / %
\setlinear \plot  5.6 -7.65 5.3 -7.55  / %
\setlinear \plot  5.6 -7.65 5.7 -7.35  / %

\setlinear \plot  3   -5.6 3.6 -4.4  / %
\setlinear \plot  3.6 -4.4 3.7 -4.7  / %
\setlinear \plot  3.6 -4.4 3.3 -4.5  / %

\setlinear \plot  4.6 -5.65 4.0 -4.4  / %
\setlinear \plot  4.6 -5.65 4.3 -5.55  / %
\setlinear \plot  4.6 -5.65 4.7 -5.35 / %

\setcoordinatesystem point at -10 0  

\put{$\circ$}[c] at 3.8 -8
\put{$+$}[c] at 7.8 -8
\put{$\bullet$}[c] at 2.8 -6
\put{$\circ$}[c] at 4.8 -6
\put{$\circ$}[c] at 6.8 -6
\put{$+$}[c] at 1.8 -4
\put{$\circ$}[c] at 5.8 -4
\put{(b) The AR-quiver of $\Gamma'$} at 4.8 -10

\setlinear \plot  4 -7.65   4.6 -6.4  / %
\setlinear \plot  4.6 -6.4 4.7 -6.7  / %
\setlinear \plot  4.6 -6.4 4.3 -6.5  / %

\setlinear \plot  3.6 -7.65 3.0 -6.4  / %
\setlinear \plot  3.6 -7.65 3.3 -7.55 / %
\setlinear \plot  3.6 -7.65 3.7 -7.35  / %

\setlinear \plot  7.6 -7.65 7.0 -6.4  / %
\setlinear \plot  7.6 -7.65 7.3 -7.55  / %
\setlinear \plot  7.6 -7.65 7.7 -7.35  / %

\setlinear \plot  5   -5.6 5.6 -4.4  / %
\setlinear \plot  5.6 -4.4 5.7 -4.7  / %
\setlinear \plot  5.6 -4.4 5.3 -4.5  / %

\setlinear \plot  2.6 -5.65 2.0 -4.4  / %
\setlinear \plot  2.6 -5.65 2.3 -5.55  / %
\setlinear \plot  2.6 -5.65 2.7 -5.35 / %

\setlinear \plot  6.6 -5.65 6.0 -4.4  / %
\setlinear \plot  6.6 -5.65 6.3 -5.55  / %
\setlinear \plot  6.6 -5.65 6.7 -5.35 / %

\endpicture
\caption{The AR-quivers of $\Gamma$ and $\Gamma'$}
\label{gammaquiver}
\end{figure}

\begin{figure}[htbp]
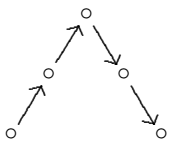

\beginpicture

\setcoordinatesystem units <0.5cm,0.4cm>             
\setplotarea x from -9 to 10, y from -8 to -1.5       

\put{$\circ$}[c] at 1.8 -8
\put{$\circ$}[c] at 5.8 -8
\put{$\circ$}[c] at 2.8 -6
\put{$\circ$}[c] at 4.8 -6
\put{$\circ$}[c] at 3.8 -4

\setlinear \plot  2 -7.65   2.6 -6.4  / %
\setlinear \plot  2.6 -6.4 2.7 -6.7  / %
\setlinear \plot  2.6 -6.4 2.3 -6.5  / %

\setlinear \plot  5.6 -7.65 5.0 -6.4  / %
\setlinear \plot  5.6 -7.65 5.3 -7.55  / %
\setlinear \plot  5.6 -7.65 5.7 -7.35  / %

\setlinear \plot  3   -5.6 3.6 -4.4  / %
\setlinear \plot  3.6 -4.4 3.7 -4.7  / %
\setlinear \plot  3.6 -4.4 3.3 -4.5  / %

\setlinear \plot  4.6 -5.65 4.0 -4.4  / %
\setlinear \plot  4.6 -5.65 4.3 -5.55  / %
\setlinear \plot  4.6 -5.65 4.7 -5.35 / %

\endpicture
\caption{The common sub-translation quiver of the AR-quivers of
$\Gamma$ and $\Gamma'$}
\label{subquiver}
\end{figure}

The AR-quivers of $\Gamma=\End_{\C}(T)^{\op}$ and $\Gamma'=\End_{\C}(T')^{\op}$
are given in Figure~\ref{gammaquiver}. The two vertices labelled
by a ``$+$'' are identified, and the simples $S_M$ and $S_{M^{\ast}}$
are shown by filled-in circles. We can see that the full sub-translation
quiver of the AR-quiver of $\Gamma$ consisting of all of the vertices except
$S_M$ is isomorphic to the full sub-translation quiver of $\Gamma'$ 
consisting of all of the vertices except $S_{M^{\ast}}$.
See Figure~\ref{subquiver}.

\begin{figure}
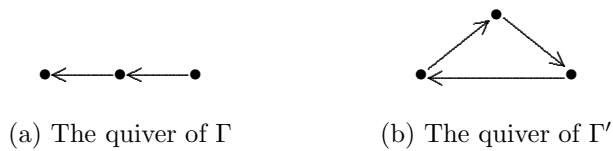

\beginpicture
\setcoordinatesystem units <0.5cm,0.4cm>             
\setplotarea x from -5 to 4, y from -3 to 4       

\put{$\bullet$}[c] at 0 0
\put{$\bullet$}[c] at 2 0
\put{$\bullet$}[c] at 4 0
\arrow <1.8mm> [.5,1] from 1.8 0 to 0.2 0
\arrow <1.8mm> [.5,1] from 3.8 0 to 2.2 0
\put{(a) The quiver of $\Gamma$} at 2 -2

\setcoordinatesystem point at -10 0  
\put{$\bullet$}[c] at 0 0
\put{$\bullet$}[c] at 2 2
\put{$\bullet$}[c] at 4 0
\arrow <1.8mm> [.5,1] from 0.2 0.2 to 1.8 1.8
\arrow <1.8mm> [.5,1] from 2.2 1.9 to 3.8 0.3
\arrow <1.8mm> [.5,1] from 3.8 -0.1 to 0.2 -0.1
\put{(b) The quiver of $\Gamma'$} at 2 -2

\endpicture
\caption{The quivers of the algebras $\Gamma$ and $\Gamma'$}
\label{a3quivers}
\end{figure}

This also gives a nice example of mutation. The quivers of
$\Gamma$ and $\Gamma'$ are shown in Figure~\ref{a3quivers};
$\Gamma$ has no relations, but for $\Gamma'$ the relations are that the
product of any pair of composable arrows is zero.
The corresponding mutation is the mutation at $2$
of the matrix $X$ to $X'$, where
$$\begin{array}{cc}
X=\left(\begin{array}{ccc}
0  & 1  & 0 \\
-1 & 0  & 1 \\
0  & -1 & 0
\end{array}\right), &
X'=\left(\begin{array}{ccc}
0  & -1  & 1 \\
1 & 0  & -1  \\
-1  & 1 & 0
\end{array}\right).
\end{array}
$$

\noindent {\bf Acknowledgements}
We would like to thank A.~Zelevinsky for many helpful discussions.
R.~Marsh would like to thank Northeastern University, Boston and NTNU,
Trondheim, for their kind hospitality.

\bibliographystyle{plain}

\end{document}